\UseRawInputEncoding
\documentclass[12pt]{amsart}
\usepackage{amssymb}
\usepackage[all]{xy}
\usepackage[toc,page,title,titletoc,header]{appendix}
\usepackage{xcolor}

\theoremstyle{plain}
\newtheorem{thm}{Theorem}[section]

\newtheorem*{thm1.1}{Theorem 1.1}

\newtheorem{lem}[thm]{Lemma}
\newtheorem{cor}[thm]{Corollary}
\newtheorem{pro}[thm]{Proposition}

\theoremstyle{definition}

\newtheorem{rem}[thm]{Remark}

\newtheorem{defi}[thm]{Definition}
\newtheorem{exe}[thm]{Example}

\newtheorem{lemdefi}[thm]{Lemma-Definition}

\numberwithin{equation}{section}
\newcounter{elno}                

\newcommand{\la}{\lambda}

\newcommand{\zar}{{\rm zar}}

\newcommand{\Supp}{{\rm Supp}}

\newcommand{\Psef}{{\rm Psef}}

\newcommand{\cone}{{\rm cone}}

\newcommand{\Hom}{{\rm Hom}}

\newcommand{\id}{{\rm id}}

\newcommand{\Eff}{{\rm Eff}}
\newcommand{\Amp}{{\rm Amp}}
\newcommand{\Bigc}{{\rm Big}}

\newcommand{\Sp}{{\rm Sp}\,}

\newcommand{\boxtensor}{{\Box\kern-9.03pt\raise1.42pt\hbox{$\times$}}}

\newcommand{\propsubset}
{\mbox{$\textstyle{
			\subseteq_{\kern-5pt\raise-1pt\hbox{\mbox{\tiny{$/$}}}}}$}}
\newcommand{\sA}{{\mathcal A}}
\newcommand{\sB}{{\mathcal B}}
\newcommand{\sC}{{\mathcal C}}
\newcommand{\sD}{{\mathcal D}}

\newcommand{\sG}{{\mathcal G}}

\newcommand{\sM}{{\mathcal M}}
\newcommand{\sN}{{\mathcal N}}
\newcommand{\sO}{{\mathcal O}}
\newcommand{\sP}{{\mathcal P}}

\newcommand{\sR}{{\mathcal R}}

\newcommand{\sZ}{{\mathcal Z}}

\newcommand{\C}{{\mathbb C}}
\newcommand{\D}{{\mathbb D}}

\renewcommand{\P}{{\mathbb P}}

\newcommand{\R}{{\mathbb R}}

\newcommand{\Z}{{\mathbb Z}}
\newcommand{\bk}{{\mathbf{k}}}

\begin{document}
	\title[]{Numerical action for endomorphisms}
	\author{Junyi Xie}
	
	\address{Beijing International Center for Mathematical Research, Peking University, Beijing 100871, China}
	
	\email{xiejunyi@bicmr.pku.edu.cn}

	\thanks{The author is supported by the NSFC Grant No.12271007.}

	\date{\today}

	\bibliographystyle{alpha}

	
	\maketitle
	

	\begin{abstract}
Let $f: X\to X$ be a surjective endomorphism of a projective variety of dimension $d$. The aim of this paper is to study the action of $f$ on the numerical group of divisors. 

For this, we introduce a notion of spectrum for an open and salient invariant cone.
Let $V$ be a finitely dimensional vector space over $\R$ and $g: V\to V$ be a linear endomorphism. Let $\sC$ be an open and salient $g$-invariant cone. We define $\Sp(g,\sC)$ to be the set of $\alpha>0$ such that no $v\in V$ satisfying $gv-\alpha v\in \sC$. We call $\Sp(g,\sC)$ the spectrum of $g$ for $\sC.$ We prove some basic properties for this spectrum. 
In particular, we get the following result: For any subset $S$ of $\C$, let $E_S(\C)\subseteq V\otimes_\R \C$ be the sum of generated eigenspaces for $g$ with eigenvalues in $S$. Then $E_S(\C)\cap \sC\neq \emptyset$ if and only if $S$ contains $\Sp(g,\sC).$ Here we view $\sC$ as a subset of $V\otimes_\R\C$ via the natural embedding $V\subseteq V\otimes_{\R}\C.$

Let $\la_i(f), i=0,\dots, d$ be the dynamical degrees and $\mu_i(f):=\frac{\la_i(f)}{\la_{i-1}(f)}, i=1,\dots,d$ be the cohomological Lyapunov exponents.  Let $\Bigc(X)$ and $\Amp(X)$ be the big and ample cones in $N^1(X)_\R.$ We show that $\Sp(f^*,\Bigc(X))=\{\mu_1(f),\dots, \mu_d(f)\}.$ We also compute the spectrum for the ample cone using the cohomological Lyapunov exponents of all periodic irreducible subvarieties. In particular, $f$ is quasi-amplified if and only if it is cohomologically hyperbolic; and $f$ is amplified if and only if every subsystem of $(X,f)$ is cohomologically hyperbolic.
As a consequence, we show that every factor of an amplified (resp. quasi-amplified) endomorphism is amplified (resp. quasi-amplified). 

We introduce a notion of generated (positive) cycles, which can be viewed as an algebraic analogy of (positively) closed currents in complex geometry. This notion plays a key role in our proof of the ample cone case. 

\end{abstract}
	
	
	\section{Introduction}

Let $\bk$ be a field. Let $X$ be a projective variety over $\bk$ of dimension $d$ and $f: X\to X$ be a surjecive endomorphism of $X$.
Denote by $N_{d-i}(X)$ the group of numerical classes of $X$ of codimension $i$ and $N^i(X):=\Hom(N_{d-i}(X),\Z)$. Set
$N_{d-i}(X)_{\R}:=N_{d-i}(X)\otimes \R$ and $N^i(X)_{\R}:=N^i(X)\otimes \R.$ Denote by $\Amp(X)$ and $\Bigc(X)$ the ample cone and the big cone in $N^1(X)_\R$. Both of them are open, $f^*$-invariant and salient i.e.  do not contain any line. The ample and big cones gives two natural notions of positivity on $N^1(X)_\R$ which contain many geometric information of $X$. 
The aim of this paper is to study the action of $f^*$ on $N^1(X)_\R.$
In the follow-up works, we will apply the results and the method developed in this paper to the Kawaguchi-Silverman conjecture and the Dynamical Mordell-Lang conjecture. 

\medskip

In general, $\dim N^1(X)_\R$ can be very large compare to $\dim X.$ A priori, the action of $f^*$ on $N^1(X)_\R$ may have many different eigenvalues. It is natural to ask, whether there is a natural  $f^*$-invariant subspace of $N^1(X)_\R$ which has smaller dimension, but still 
captures a large part of dynamical information.  For this purpose, we want such space cuts the ample (resp. big) cone. For this, we introduce a notion of spectrum for a good cone. 

\subsection{Spectrum for a good invariant cone}
Let $W$ be a finite dimensional $\R$-vector space and $g:W\to W$ be an endomorphism.  
Set $W_{\C}:=W\otimes_{\R}\C$ and $g_{\C}: W_{\C}\to W_{\C}$ the endomorphism induced by $g.$
View $W$ as a $\R$-subspace of $W_{\C}$. We have $g=g_{\C}|_W$.
Denote by $\Sp(g)$ the set of eigenvalues of $g.$


\subsubsection{Generated eigenspaces}	For every $c\in \C$, denote by $E_c(\C)$ the generated eigenspace of $g_{\C}: W_{\C}\to W_{\C}$ i.e. $$E_c(\C):=\cup_{n\geq 0}\ker((g-c\id)^n:W_{\C}\to W_{\C}).$$ 
It is clear that $E_c(\C)\neq 0$ if and only if $c\in \Sp(g).$
When $c\in \R$, write $E_c:=E_c(\C)\cap W.$ We have $E_c(\C)=E_c\otimes_{\R}\C.$
For every $S\subseteq \C$, set $$E_S(\C):=\oplus_{c\in S}E_c \text{ and } E_S:=E_S(\C)\cap W.$$
We note that, if $S=\overline{S}$, then $E_S(\C)=E_S\otimes_{\R} \C.$ 
We mainly interest in the case where $S\subseteq \R.$ In this case $E_S=\oplus_{c\in S}E_c.$

\medskip

\subsubsection{Good invariant cone}
Let $\sC$ be a non-empty open convex cone in $W$ satisfying the following properties:
\begin{points}
	\item $g(\sC)=\sC$;
	\item $\overline{\sC}$ is salient i.e. $\overline{\sC}\setminus \{0\}$ is convex.
\end{points}
We call such $\sC$ a \emph{good invariant cone} for $g.$

\begin{exe} If $W=N^1(X)_{\R}$ for a projective variety $X$, $g=f^*: N^1(X)_{\R}\to N^1(X)_{\R}$ for some surjective endomorphism $f: X\to X$, then the ample cone $\Amp(X)$ and the big cone $\Bigc(X)$ are good invariant cone for $g$.
\end{exe}

Let $\sC$ be a good invariant cone for $g$.

\begin{defi}For $\alpha\in \R$, we say that $g$ is \emph{$\alpha$-amplified} for $\sC$, if there is $N\in W$ such that $gN-\alpha N\in \sC$.
	Define the $\sC$-spectrum $\Sp(g,\sC)$ for $g$ to be the set of $\alpha\in \R$ such that $g$ is {\bf not} $\alpha$-amplified.
\end{defi}

\medskip

If $\alpha\not\in \Sp(g)$, then $g-\alpha$ is invertible on $W$, hence $g$ is $\alpha$-amplified for $\sC$. In particular, we have $ \Sp(g,\sC)\subseteq \Sp(g).$ It is clear that $\Sp(g,\sC)$ is decreasing on $\sC$.
The following result gives a description of the spectrum $\Sp(g, \sC)$ using generated eigenspaces.
\begin{thm}\label{thmceighearintro}
	For every subset $S\subseteq \C$, $\sC\cap E_S\neq\emptyset$ if and only if $\Sp(g,\sC)\subset S.$ 
\end{thm}

\medskip

\subsection{Cohomological Lyapunov exponents}\label{subsecdynamicaldeg}
For every $i=0,\dots,d$, the \emph{$i$-th dynamical degree} $\la_i(f)$ of $f$ is 
is defined to be the spectral radius of 
$$f^*: N^i(X)_\R\to N^i(X)_\R.$$
It can be also defined as follows: Let $L$ be any ample line bundle on $X$, then $$\la_i(f)=\lim_{n\to \infty} ((f^n)^*L^i\cdot L^{d-i})^{1/n}.$$
The later definition can be generated to rational self-maps \cite{Russakovskii1997,Dinh2005,Dinh2004,Truong2020,Dang2020}.

By \cite[Theorem 1.1(3)]{Truong2020} (see also \cite{Dinh2005,Dang2020}), the sequence $\la_i(f), i=0,\dots, d$ is log-concave i.e. 
$$\la_i(f)^2\geq \la_{i-1}(f)\la_{i+1}(f)$$ for $i=1,\dots, d-1.$
As in \cite{Xie2024}, we define $\mu_i(f):=\la_i(f)/\la_{i-1}(f)$ for $i=1,\dots, d$ and $\mu_{d+1}(f):=0.$ We call them the cohomological Lyapunov exponents.
The log-concavity of $\la_i(f), i=0,\dots, d$ implies that the sequence $\mu_i(f), i=1,\dots, d+1$  is decreasing.

\subsection{Main results}
%
%
%
%
For $\alpha\in \R_{>0}$, we say that $f$ is \emph{$\alpha$-quasi-amplified} (resp. \emph{$\alpha$-amplified}) if $$\alpha\not\in \Sp(f^*,\Bigc(X) ) \,\,(\text{resp.  }    \alpha\not\in \Sp(f^*,\Amp(X)))$$ i.e. there is $N\in N^1(X)_{\R}$ such that $f^*N-\alpha N$ is big (resp. ample). 

\medskip

As defined in \cite[Definition 2.2]{Meng2023a} and in \cite{fa},  $f$ is called \emph{quasi-amplified} (resp. amplified), if there is $N\in N^1(X)_{\R}$ such that $f^*N-N$ is big (resp. ample). So $f$ is quasi-amplified (resp. amplified) if and if $f$ is $1$-quasi-amplified (resp. $1$-amplified).

\medskip

The following result computes the spectrum for the big cone. 
\begin{thm}\label{thmchendamplifyintro}
	We have $$\Sp(f^*,\Bigc(X))=\{\mu_i(f) |\,\, i=1,\dots, d\}.$$
	In other words, for $\alpha\in \R_{>0}$, $f$ is $\alpha$-quasi-amplified if and only if $$\alpha\not\in\{\mu_i(f) |\,\, i=1,\dots, d\}.$$
	In particular, $f$ is quasi-amplified if and only if $f$ is cohomologically hyperbolic. 
\end{thm}

As a consequence of Theorem \ref{thmchendamplifyintro}, every $\mu_i(f), i=1,\dots, d$ is an eigenvalue of $f^*: N^1(X)_\R\to N^1(X)_\R.$

The proof of Theorem \ref{thmchendamplifyintro} relies on the recursive inequalities constructed in \cite[Theorem 3.7]{Xie2024} and the computation of mixed degrees in \cite[Corollary 3.4]{Xie2024}.

\medskip

Next, we compute the spectrum for the ample cone. 
Let $V$ be an irreducible and periodic subvariety of $X$ of dimension $d_V\geq 0$. Define $$\mu_i(V,f):=\mu_i(f|_V^{r_V}:V\to V)^{1/r_V}$$ where $r_V\geq 1$ is a period of $V$. It does not depend on the choice of $r_V.$

\begin{thm}\label{thmquaisamptoamintr}
	We have $$\Sp(f^*,\Amp(X))=\cup_{V}\{\mu_i(V,f) |\,\, i=1,\dots, d_V\}$$ where the union taken over all irreducible periodic subvarieties.

	In other words,	for $\alpha\in \R_{>0}$, $f$ is $\alpha$-amplified if and only if for every periodic irreducible subvariety $V$, $f^{r_V}|_V$ is $\alpha$-quasi-amplified, where $r_V\geq 1$ is a period of $V$. 
\end{thm}

To prove Theorem \ref{thmquaisamptoamintr}, we introduce a notion of \emph{generated (positive) cycles} (c.f. Section \ref{sectiongencycle}), which can be viewed as the algebraic analogy of (positively) closed currents in complex geometry. 
We prove some basic properties of it. In particular, we prove a decomposition theorem of generated positive cycles (c.f. Theorem \ref{thmdecogc}), which can be viewed as an algebraic analogy of Siu's decommposition theorem \cite{Siu1974}.

\medskip

As a consequence of Theorem \ref{thmchendamplifyintro} and Theorem \ref{thmquaisamptoamintr}, we get the following result.
\begin{cor}[c.f. Corollary \ref{corfactoramplified} and \ref{corfactoramplifiedquasi}]\label{corfactoramplifiedboth}Let $Y$ be a projective variety over $\bk$ and $g:Y\to Y$ be an endomorphism. Let $\pi: X\to Y$ be a surjective morphism such that $\pi\circ f=g\circ \pi.$ If $f$ is $\alpha$-amplified (resp. $\alpha$-quasi-amplified) for some $\alpha\in \R_{>0}$, then $g$ is $\alpha$-amplified (resp. $\alpha$-quasi-amplified).
\end{cor}

\subsubsection{Hyperbolicity}
The case $\alpha=1$ in Theorem \ref{thmchendamplifyintro} and \ref{thmquaisamptoamintr} are especially interesting as they clarify two different algebraic analogies of the notion of hyperbolicity.

\medskip

The hyperbolicity is one of the most important notion in smooth dynamical system. Basically, it means that the tangent bundle can be decomposed to the direct sum of the attracting and the repelling parts. Further, if we only have the repelling part, such map is called \emph{expanding}.
These notion are analytic, which do not make sense for algebraic dynamical system over an abstract field. Indeed, even when $\bk=\C$, 
very few algebraic dynamical system could be Anosov (which is a strong version of hyperbolicity) c.f. \cite{Ghys1995,Cantat2004,Xu2024}.

\medskip

In our setting (or more generally for rational self-maps), $f$ is called \emph{cohomologically hyperbolic} if the cohomological Lyapunov exponents $\mu_i(f)\neq 1$ for every $i=1,\dots, d.$ From the point of view of ergodic theory, this notion gives an algebraic analogy of the hyperbolicity.

\medskip

In algebraic geometry, we may view an ample line bundle as an analogy of a Riemannian metric.
From this point of view, the algebraic analogy of an expending map should be an \emph{int-amplified endomorphism}\footnote{An endomorphism $f: X\to X$ is called int-amplified if there is an ample line bundle $L$ of $X$ such that $f^*L-L$ is ample \cite{Meng2020}.}. It was observed by Matsuzawa, that $f$ is int-amplified if and only if $\mu_d(f)>1$ (c.f. \cite[Proposition 3.7]{Meng2023b}).
In other words, the metric and ergodic theory styles analogies of expending maps are the same.

Further, we view an arbitrary line bundle as an analogy of a pseudo-Riemannian metric.  From this point of view, the algebraic analogy of a hyperbolic map should be an amplified endomorphism. 
Then the $\alpha=1$ case of Theorem \ref{thmchendamplifyintro} and 
\ref{thmquaisamptoamintr} connects the metric and ergodic theory styles analogies of the hyperbolicity.

\subsection*{Acknowledgement}
The author would like to thank Sheng Meng and De-Qi Zhang for their helpful comments for the first version of the paper.

	\section{Linear algebra for good invariant cones}\label{sectionlinearalg}
In this section $W$ is a finite dimensional $\R$-vector space and $g:W\to W$ be an endomorphism.  
Set $W_{\C}:=W\otimes_{\R}\C$ and $g_{\C}: W_{\C}\to W_{\C}$ the endomorphism induced by $g.$
View $W$ as a $\R$-subspace of $W_{\C}$. We have $g=g_{\C}|_W$.
Denote by $\Sp(g)$ the set of eigenvalues of $g.$ 
	
\subsection{Speed of growth}
%
%
	Let $\|\cdot\|$ by any norm on $W$. Then for every $v\in W$, we have 
	\begin{flalign}\label{equnormgrowth}
		\lim_{n\to \infty}\|f^n(v)\|^{1/n}=&\min\{r\in \R_{\geq 0}|\,\, v\in E_{\overline{\D}_{r}}\}\\
		\label{equnspgrothrate}=&\min\{|c||\,\, c\in \Sp(g), v\in E_{\overline{\D}_{|c|}}\}.
		\end{flalign}
Here $\overline{\D}_{r}$ is the closed disc in $\C$ with center $0$ and radius $r.$
	
	\medskip
	
The following lemma is useful. It's proof is a simple application of the Jordan normal form. We leave it to the readers.
	\begin{lemdefi}\label{lemgrowthspeed}
Assume that $\Sp(g)\subseteq \R_{>0}.$ Then for every $Z\in W^{\vee}$ and $v\in W$, there is a unique $(\beta(v),a(v))\in \R_{>0}\times \Z_{\geq 0}$ such that $(g^n(v)\cdot Z)=C(v)\beta(v)^nn^{a(v)}+O(\beta(v)^nn^{a(v)-1})$ with $C(v)\neq 0.$
Moreover, there is a unique $(\beta,a)\in \R_{>0}\times \Z_{\geq 0}$ such that the following holds:
\begin{points}
	\item for every $v\in W$, $(\beta,a)\geq (\beta(v),a(v))$ for the lexicographical order $\geq$;
	\item $\{v\in W|\,\, (\beta,a)>(\beta(v),a(v))\}$ is a proper closed subspace of $W.$
\end{points}
We say that $v\in W$ has \emph{maximal growth} for $Z$ if $(\beta,a)=(\beta(v),a(v))$.
\end{lemdefi}

Let $\sC$ be a non-empty open convex cone in $W$ satisfying the following properties:
\begin{points}
	\item $g(\sC)=\sC$;
	\item $\overline{\sC}$ is salient i.e. $\overline{\sC}\setminus \{0\}$ is convex.
\end{points}
We call such $\sC$ a \emph{good invariant cone} for $g.$

\begin{lem}\label{lemgoodconegr}
	Let $Z\in W^{\vee}$. Assume that for every $v\in \overline{\sC}\setminus \{0\}$, $(v, Z)>0.$
	Then $\log \frac{(g^n(v),Z)}{\|g^n(v)\|}$ is bounded.
	In particular, if $\Sp(g)\subseteq \R_{>0},$ then every $v\in \sC$ has maximal growth for $Z$.  
	\end{lem}
\proof
The cone $\sC$ with $Z$ induces a norm $\|\cdot\|'$ on $W$ as follows: for every $v\in W$, $\|v\|':=\inf(Z(v_1)+Z(v_2))$ where $(v_1,v_2)$
are taken over all pairs in $\overline{\sC}$ with $v_1-v_2=v.$ Easy to check that $\|\cdot\|$ is a norm and for every $v\in \overline{\sC}$, $\|v\|'=Z(v).$ As all norms on $W$ are equivalent, we conclude the proof.
\endproof

\subsection{Spectrum for good invariant cones}
%
%
%

%

\begin{lem}\label{lemrealspace}
Let $\sC$ be a convex subset of $W$ which spans $W$.
Let $V$ be a $g$-invariant subspace of $W.$ Then 
$V\cap \sC\neq \emptyset$ if and only if $V\cap E_{\R_{>0}}\cap \sC\neq \emptyset$.
\end{lem}

\proof
The ``if" parts are trivial. We now prove the ``only if" parts. 
After replacing $W$  by the subspace spans by $\sC$, we may assume that $\sC$ spans $W$.
After replacing $\sC$ by $\cone(\sC):=\R_{>0}\sC$, we may assume that $\sC^{\circ}\neq \emptyset.$
As $\sC^{\circ}\neq \emptyset$ and $g(\sC)=\sC$, $g$ in invertible i.e. $E_0=0.$

\medskip

Pick $L\in V\cap \sC$.
By contradiction, assume that $V\cap E_{\R}\cap \sC=\emptyset$. 
Then $$V\cap E_{\R}\cap \cone(\sC)=\emptyset.$$
By Hahn-Banach theorem, there is $Z\in W^{\vee}$ such that 
$V\cap E_{\R_{>0}}\subseteq Z^{\bot}$ and for every $N\in \sC$, $(N\cdot Z)> 0$.

\medskip

Set $S:=\Sp(f^*|_V).$ We have $S=\overline{S}.$
Consider the canonical decomposition in $V\otimes_{\R}\C$ that $$L=\sum_{c\in S} L_c$$ where $L_c\in E_c(\C).$
It is clear that $L_{\overline{c}}=\overline{L_{c}}.$ For $n\geq 0$, write $L_n:=g^n(L)$ and $(L_c)_n:=g^n(L_c)$.
As $V\cap E_{\R}\subseteq Z^{\bot}$, we have
$$(L_n\cdot Z)=\sum_{c\in S\setminus \R_{>0}}((L_c)_n\cdot Z).$$
As $L\in V\cap \sC$ and for every $N\in \sC$, $(N\cdot Z)> 0$, we have 
\begin{equation}\label{equationhpositive}h(n):=(L_n\cdot Z)> 0\end{equation} for every $n\geq 0.$
Easy to see that $h(n)$ is an exponential-polynomial function i.e. $h(n)$ is a finite sum of terms having form 
$u^nn^s$
where $u\in \C^*$ and $s\in \Z_{\geq 0}.$ Moreover, we may ask the $u$ above are contained in $\Sp(g)\setminus\R_{>0}.$ We note that $0\not\in \Sp(g)\setminus\R_{>0}.$
Consider $I:=(\Sp(g)\setminus \R_{>0})\times \Z_{\geq 0}$. 
For every $p=(u,s)\in I$, set $H_p(n):=u^nn^s.$
Then there is a finite set $S\subseteq I$ such that
$$h(n)=\sum_{p\in S}c_pH_p(n)$$ where $c_p\in \C^*.$
By (\ref{equationhpositive}), $h$ is not constantly equal to $0$, hence $S\neq \emptyset.$
As $h(n)\in \R$ for every $n\geq 0$, $S$ is invariant under the complex conjugacy i.e. for every $p=(u,s)\in S$, we have $\overline{p}:=(\overline{u},s)\in S$. Moreover, we have $c_{\overline{p}}=\overline{c_p}.$

\medskip

Set $J:=\R_{>0}\times \Z_{\geq 0}$ with the lexicographical order $\geq$. For $p=(u,s)\in I$, write $|p|:=(|u|,s)\in J.$
There is a unique maximal element $(\beta,a)\in \{|p||\,\, p\in S\}.$ 
Set $S^+:=\{p\in S|\,\, |p|=(\beta,a)\}$ and $S^-=S\setminus S^+.$ Then $S^+, S^-$ are invariant under complex conjugacy.
Set $h^+:=\sum_{p\in S^+}c_pH_p$ and $h^-:=\sum_{p\in S^-}c_pH_p.$ We have $h=h^++h^-.$
There is $C>0$ such that 
\begin{equation}\label{equhmestihm}
	|h^-(n)|\leq C'\beta^nn^{a-1}.
\end{equation}
for every $n\geq 1.$

For every $p\in S^+$, write $$H_p(n)=e^{i\theta_pn}\beta^nn^a$$ where $\theta_p\in \R/2\pi \Z.$
Note that, we have $\theta_p\neq 0$ for every $p\in S^+.$
The map $$p\in S^+\mapsto \theta_p\in(\R/2\pi \Z)\setminus \{0\}$$ is injective. 
Set $$C(n):=\sum_{p\in S^+}c_pe^{i\theta_pn}.$$
For every $n\geq 0$, we have $C(n)\in \R.$
We have  $h^+(n)=C(n)\beta^nn^a.$ Using Vandermonde determinant, easy to show that $C(n)$ is not constantly equal to $0.$

\medskip

By Poincar\'e recurrence theorem,  there is a strict increasing sequence $n_i\in \Z_{\geq 0}, i\geq 0$ such that 
$$\lim_{i\to \infty}\theta_pn_i= 0\in \R/2\pi \Z$$ for every $p\in S^+.$

\begin{lem}\label{lemclneg} There is $m\geq 0$ such that $C(m)<0.$
	\end{lem}
Set $m_i:=m+n_i,$ then $m_i$ is strictly increasing and $C(m_i)\to C(m)$ as $i\to \infty.$
By (\ref{equhmestihm}),
we get $$\lim_{i\to\infty}\frac{h(m_i)}{\beta^{m_i}m_i^a}=\lim_{i\to \infty}\frac{h^+(m_i)}{\beta^{m_i}m_i^a}=\lim_{i\to \infty}C(m_i)= C(m)<0.$$
This contradicts with (\ref{equationhpositive}).
\endproof

\proof[Proof of Lemma \ref{lemclneg}]
By contradition, we assume that $C(n)\geq 0$ for every $n\geq 0.$
There is $m'\geq 0$ such that $C(m')\neq 0.$
So we have $C(m')>0.$

\medskip

For every $l\geq 0$, we have 
\begin{equation}\label{equationsumcnup}
	0\leq \sum_{n= 0}^lC(n)=\sum_{p\in S^+}c_p\sum_{n=0}^le^{i\theta_pn}=\sum_{p\in S^+}c_p\frac{1-e^{i\theta_p(l+1)}}{1-e^{i\theta_p}}\leq\sum_{p\in S^+}\frac{2|c_p|}{1-e^{i\theta_p}}.
\end{equation}

As $$\lim_{i\to \infty}\theta_pn_i= 0\in \R/2\pi \Z,$$ for every $p\in S^+,$
$$C(m'+n_i)\to C(m)$$ as $i\to \infty.$
Then we have 
$$\liminf_{j\to \infty}\sum_{n= 0}^{m'+n_j}C(n)\geq \lim\sum_{i=0}^{j}C(m'+n_i)=+\infty.$$
This contradicts with (\ref{equationsumcnup}). We concludes the proof.
\endproof

Let $\sC$ be a good invariant cone for $g$.

\begin{defi}For $\alpha\in \R$, we say that $g$ is \emph{$\alpha$-amplified} for $\sC$, if there is $N\in W$ such that $gN-\alpha N\in \sC$.
Define the $\sC$-spectrum $\Sp(g,\sC)$ for $g$ to be the set of $\alpha\in \R$ such that $g$ is {\bf not} $\alpha$-amplified.
\end{defi}

\medskip

If $\alpha\not\in \Sp(g)$, then $g-\alpha$ is invertible on $W$, hence $g$ is $\alpha$-amplified for $\sC$. In particular, we have $ \Sp(g,\sC)\subseteq \Sp(g).$ It is clear that $\Sp(g,\sC)$ is decreasing on $\sC$.

\medskip

\begin{lem}\label{lemquasidyna}For $\alpha\in \R$, the following statement are equivalent:
	\begin{points}
		\item $g$ is $\alpha$-amplified for $\sC$;
		\item there is $n\geq 1$, such that $g^n$ is $\alpha^n$-amplified for $\sC$;
		\item for every $n\geq 1$, such that $g^n$ is $\alpha^n$-amplified for $\sC$.
	\end{points}
	\end{lem}
\proof
It is clear that (iii) implies (i) and (i) implies (ii).

\medskip
Now we show that (i) implies (iii).
By (i), there is $N\in W$ such that $$L:=g(N)-\alpha N\in \sC.$$ 
For every $n\geq 1$, we have 
$$g^n(N)-\alpha^n N=\sum_{j=1}^{n}(\alpha^{n-j}g^j(N)-\alpha^{n-j+1} g^{j-1}(N))=\sum_{j=0}^{n-1}g^j(L)\in \sC.$$
This implies (iii).

\medskip

We only need to show that (ii) implies (i). 
By (ii), there is $n\geq 1$ and $N\in W$ such that $$M:=g^n(N)-\alpha^n N\in \sC.$$ 
Set $N':=\sum_{j=0}^{n-1}\alpha^{n-j}g^j(N).$ Then 
$$f^*N'-\alpha N'=\alpha g^n(N)-\alpha^{n+1}N=\alpha M\in \sC.$$ 
This implies (i).
\endproof

Lemma \ref{lemquasidyna} has the following direct consequence.
\begin{cor}\label{corspgnandg}	For every $n\geq 1$, $\Sp(g^n,\sC)=\{\alpha^n|\,\, \alpha\in \Sp(g,\sC)\}.$
	\end{cor}

The following lemma shows that $\Sp(g,\sC)\subseteq \R_{>0}\cap \Sp(g)$ and gives constrains of $g$-invariant subspace $V$ meeting $\sC.$
We will reinforce this lemma in Corollary \ref{corinnersubsamesp} latter.

\begin{lem}\label{leminvsp}
Let $V$ be a $g$-invariant subspace such that $V\cap \sC\neq\emptyset$, then we have $\Sp(g,\sC)\subseteq\Sp(g|_{V},\sC\cap V).$
As a consequence, we have $\Sp(g,\sC)\subseteq \R_{>0}$ and $\dim V\geq \#\Sp(g,\sC).$
\end{lem}
\proof
By contradiction, assume that $\sC\cap V\neq\emptyset$ and $\Sp(g,\sC)\not\subseteq \Sp(g|_{V}).$ Pick $c\in S\setminus \Sp(g,\sC).$ As $g-c\id$ is invertible on $V$, $V\subseteq {\rm Im}(g-c\id)$. So ${\rm Im}(g-c\id)\cap \sC$ contains $\sC\cap V\neq\emptyset$.
By Lemma \ref{lemrealspace}, we have $\Sp(g,\sC)\subseteq \R_{>0}.$
\endproof

We give an example that $\Sp(g,\sC)=\R_{>0}\cap \Sp(g).$

\begin{exe}
	Let $W=\R^d$ with standard base $e_1,\dots,e_d$. Let $a_1,\dots, a_d\in \R_{>0}.$ Let $g: W\to W$ be the morphisms sending $e_i$ to $a_ie_i$, $i=1,\dots, d.$ Let $$\sC:=\{x_1e_1+\dots+x_de_d|\,\, x_1,\dots,x_d\in \R_{>0}\}.$$ Then $\sC$ is a good invariant cone for $g$. Easy to check that $\Sp(g,\sC)=\{a_1,\dots, a_d\}=\R_{>0}\cap \Sp(g).$
\end{exe}

The following result gives a description of the spectrum $\Sp(g, \sC)$ using generated eigenspaces.

\begin{thm}(=Theorem \ref{thmceighearintro})\label{thmceighear}
For every subset $S\subseteq \C$, $\sC\cap E_S\neq\emptyset$ if and only if $\Sp(g,\sC)\subset S.$ 
\end{thm}

\proof
The ``only if" part follows from Lemma \ref{leminvsp}.
We now prove the ``if" part. 
We only need to show that when $S=\Sp(g,\sC)$, $\sC\cap E_S=\emptyset.$

By contradiction, assume that $\sC\cap E_S=\emptyset.$ 
By Hahn-Banach theorem, there is $Z\in W^{\vee}$ such that 
$E_S\subseteq Z^{\bot}$ and for every $N\in \sC$, $(N\cdot Z)> 0$.
By Lemma \ref{lemrealspace}, $\sC\cap E_{\R_{>0}}\neq\emptyset.$
As $\Sp(g|_{E_{\R_{>0}}})\subseteq \R_{>0}$, by Lemma \ref{lemgoodconegr}, every $L'\in \sC\cap E_{\R_{>0}}$ has maximal growth for $Z.$ More precisely, there is $(\beta,a)\in \R_{>0}\times \Z_{\geq 0}$ and $C'\in \R^*$ such that $$(L_n'\cdot Z)=C'\beta^nn^a+O(\beta^nn^{a-1}),$$
where $L_n':=g^n(L').$
Moreover, for every $N\in W$, $(g^n(N)\cdot Z)=O((L_n'\cdot Z)).$

After replacing $g$ by $\beta^{-1}g$, we may assume that $\beta=1.$
As $E_S\subseteq Z^{\bot}$, $1\not\in S=\Sp(g,\sC)$. So there is $N\in W$ such that 
$$g(N)-N=L$$ for some $L\in \sC.$
Set $N_n:=g^n(N), L_n:=g^n(L)$ for $n\geq 0.$
Then  $$(L_n,Z)=C\beta^nn^a+O(\beta^nn^{a-1}),$$ for some $C\in \R^*.$
We have 
\begin{equation}\label{equationndiffel}N_{n+1}-N_n=L_n.
	\end{equation}
As $L$ has maximal growth for $Z,$
there is $m\geq 0$ such that for every $n\geq m$,
\begin{equation}\label{equationhngrowth}C/2<\frac{h(n)}{n^a}<2C.
\end{equation}
As $L_n\in \sC$, $h(n)>0.$ Hence $C>0.$
By (\ref{equationndiffel}), for every $n\geq m$, we get 
\begin{equation}\sum_{i=m}^nL_i=N_{n+1}-N_m
	\end{equation}
There is $B>0$ such that $BL\pm N\in \sC.$ Then we get 
$$\sum_{i=m}^nh(i)=\sum_{i=m}^n(L_i\cdot Z)\leq B((L_{n+1}\cdot Z)+(L_m\cdot Z))\leq B(h(n+1)+h(m)).$$
By (\ref{equationhngrowth}), we get 
$$(C/2)\sum_{i=m}^ni^a\leq 2BC((n+1)^a+m^a).$$
Let $n\to \infty$, we get a contradiction. This concludes the proof.
\endproof

\begin{cor}\label{corinnersubsamesp}
	Let $V$ be a $g$-invariant subspace such that $V\cap \sC\neq\emptyset$, then we have $\Sp(g,\sC)=\Sp(g|_{V},\sC\cap V).$
	\end{cor}
\proof
By Lemma \ref{leminvsp}, we have $\Sp(g,\sC)\subseteq \Sp(g|_{V},\sC\cap V).$
We only need to show the inverse direction.
Apply Theorem \ref{thmceighear} to $g|_V$, we get 
$$E_{\Sp(g|_{V},\sC\cap V)}\cap V\cap \sC=(E_{\Sp(g|_{V},\sC\cap V)}\cap V)\cap (V\cap \sC)\neq\emptyset.$$
Hence $E_{\Sp(g|_{V},\sC\cap V)}\cap \sC\neq\emptyset.$
We conclude the proof by Theorem \ref{thmceighear}.
\endproof


\section{The spectrum for the big cone}
Let $X$ be a projective variety over $\bk$ of dimension $d$ and $f: X\to X$ is a surjective endomorphism. 
Denote by $\mu_i:=\mu_i(f), i=1,\dots, d+1$, the cohomological Lyapunov exponents of $f.$

\medskip

%
%
%

\begin{thm}(=Theorem \ref{thmchendamplifyintro})\label{thmchendamplify}
	We have $$\Sp(f^*,\Bigc(X))=\{\mu_i |\,\, i=1,\dots, d\}.$$
	In other words, for $\alpha\in \R_{>0}$, $f$ is $\alpha$-quasi-amplified if and only if $$\alpha\not\in\{\mu_i |\,\, i=1,\dots, d\}.$$
	In particular, $f$ is quasi-amplified if and only if $f$ is cohomologically hyperbolic. 
\end{thm}

\proof[Proof of Theorem \ref{thmchendamplify}]
Let $L$ be an ample line bundle on $X$. Set $L_n:=(f^n)^*L.$

\medskip

We first assume that $\alpha\not\in\{\mu_i |\,\, i=1,\dots, d\}$.
There is a minimal $i=0,\dots, d$ such that $\mu_{i+1}<\alpha.$
If $i=0$, then $\alpha>\mu_1.$ As $\mu_1$ is the spectral radius of $$f^*: N^1(X)_{\R}\to N^1(X)_{\R},$$ $f^*-\alpha\id$ is invertible on $N^1(X)_{\R}.$
So for every big class $L\in N^1(X)_{\R}$, there is $N\in N^1(X)_{\R}$ such that $$f^*N-\alpha N=L.$$
Now assume that $i\geq 1.$ Then we have $\alpha\in (\mu_i, \mu_{i+1}).$
Pick $\epsilon\in (\max\{\alpha/\mu_i, \mu_{i+1}/\alpha\},1).$
There is $m_0\geq 1$, such that  for every $m\geq m_0,$ 
\begin{equation}\label{inequmuialpeps}
	\epsilon^m\mu_i^{m}-\alpha^m-\mu_i^{m}\mu_{i+1}^m\alpha^{-m}>0.
\end{equation}
By \cite[Theorem 3.7]{Xie2024},  there is $m_{1}>m_0$, such that for every $m\geq m_{1},$
$$M:=L_{2m}+\mu_i^{m}\mu_{i+1}^mL-\epsilon^m\mu_i^{m}L_m$$ is big.
Set $N:=L_m-\mu_i^m\mu_{i+1}^m\alpha^{-m}L.$
Then we have 
\begin{flalign}
	(f^m)^*N-\alpha^mN&=(L_{2m}-\alpha^m L_m)-\mu_i^{m}\mu_{i+1}^m\alpha^{-m}(L_m-\alpha^m L)\\
	&=M+(\epsilon^m\mu_i^{m}-\alpha^m-\mu_i^{m}\mu_{i+1}^m\alpha^{-m})L_m.
\end{flalign}
By (\ref{inequmuialpeps}), $(f^m)^*N-\alpha^mN$ is big. So $f^m$ is $\alpha^m$-quasi-amplified. By Lemma \ref{lemquasidyna}, $f$ is $\alpha$-quasi-amplifed.

\medskip

Now assume that $\alpha=\mu_i$ for some $i=1,\dots, d$ and want to show that $f$ is not $\alpha$-quasi-amplifed. 
Otherwise, we assume that $f$ is $\alpha$-quasi-amplifed, then there is $N\in N^1(X)_{\R}$ such that $$M:=f^*N-\alpha N=f^*N-\mu_i N$$ is big. 
As $M$ is big, after replacing $N$ by a suitable multiple, we may assume that $M\geq L$ i.e. $M-L$ is peudo-effective.
There is $B>0$ such that 
$$-BL\leq N\leq BL.$$
For $n\geq 0$, set $N_n:=(f^n)^*N$ and $M_n:=(f^n)^*M.$
Set $$h(n,m):=\mu_i^{-m}(L_n^{i-1}\cdot L_m\cdot L^{d-i}).$$
For $m,n\in \Z_{\geq 0}$ and $m_0\geq 0$, we have
\begin{equation}\label{equalineqsum}\sum_{j=m_0}^m\mu_i^{-j}(L_n^{i-1}\cdot M_j\cdot L^{d-i})=\mu_i^{-m-1}(L_n^{i-1}\cdot N_{m+1}\cdot L^{d-i})-\mu_i^{-m_0}(L_n^{i-1}\cdot N_{m_0}\cdot L^{d-i}).
\end{equation}
As $M\geq L$, we have 
\begin{equation}\label{equalhssum}
	\sum_{j=m_0}^m\mu_i^{-j}(L_n^{i-1}\cdot M_j\cdot L^{d-i})\geq \sum_{j=m_0}^m h(n,j).
\end{equation}
As $-BL\leq N\leq BL,$ we have 
\begin{equation}\label{equarhsdif}
	\mu_i^{-m-1}(L_n^{i-1}\cdot N_{m+1}\cdot L^{d-i})-\mu_i^{-m_0}(L_n^{i-1}\cdot N_{m_0}\cdot L^{d-i})\leq Bh(n,m+1)+B h(n,m_0)
\end{equation}
Combining (\ref{equalineqsum}), (\ref{equalhssum}) and (\ref{equarhsdif}),
we have 
\begin{equation}\label{equhineqwilcon}\sum_{j=m_0}^m h(n,j)\leq Bh(n,m+1)+B h(n,m_0).
\end{equation}
One may check that $h:\Z\geq 0\times\Z_{\geq 0}\to \R\subseteq \C$ is an exponential-polynomial function.
\begin{lem}\label{lemexppoly}Let $h:\Z\geq 0\times\Z_{\geq 0}\to \C$ be an exponential-polynomial function i.e. $h(n,m)$ is a finite sum of terms having form 
	$$u^nv^mn^sm^t$$
	where $u,v\in \C^*$ and $s,t\in \Z_{\geq 0}.$
	Assume that $h$ has the following conditions:
	\begin{points}
		\item $h$ is real and positive i.e. $h(n,m)\in \R_{>0}$ for every $n,m\in \Z_{\geq 0};$
		\item there is $D>0$ such that $$\max\{h(n+1,m),h(n,m+1)\}\leq Dh(n,m)$$ for every $n,m\in \Z_{\geq 0}.$
	\end{points}
	Then there is $C>1$, $\epsilon_0\in(0,1)$, $\beta,\gamma>0$, $a,b\in \Z_{\geq 0}$ such that 
	$$C^{-1}\leq \frac{h(n,m)}{\beta^n\gamma^mn^am^b}\leq C$$
	for $$(n,m)\in \sN:=\{(n,m)\in \Z_{\geq 0}|\,\,  (\log n)^2 \leq m\leq n^{\epsilon_0}\}.$$
	Moreover, if the following holds,
	\begin{points}
		\item[\rm(iii)] there is $\la>0$ such that for every $\delta\in (0,1)$, there is $D_{\delta}>1$ such that for every $n\geq m\geq 0$, we have 
		$$D_{\delta}^{-1}\delta^n\leq \frac{h(n,m)}{\la^n}\leq D_{\delta}\delta^{-n},$$
	\end{points}
	then we may take $\beta=\la$ and $\gamma=1.$
\end{lem}
It is clear that assumptions (i) in Lemma \ref{lemexppoly} are satisfied for our $h$.
As $L$ is ample, there is $D'>0$ such that $f^*L\leq D'L.$
So we have 
\begin{flalign}\label{equiiione}
	h(n+1,m)&=\mu_i^{-m}(L_{n+1}^{i-1}\cdot L_m\cdot L^{d-i})\\
	&= \mu_i^{-m}((f^n)^*(f^*L)^{i-1}\cdot L_m\cdot L^{d-i})\\
	&\leq D'^{i-1}\mu_i^{-m}((f^n)^*L^{i-1}\cdot L_m\cdot L^{d-i})\\
	&=D'^{i-1}h(n,m).
\end{flalign}
and 
\begin{flalign}\label{equiiione}
	h(n,m+1)&=\mu_i^{-m-1}(L_n^{i-1}\cdot L_{m+1}\cdot L^{d-i})\\
	&= \mu_i^{-m-1}(L_n^{i-1}\cdot (f^m)^*(f^*L)\cdot L^{d-i})\\
	&\leq (D'/\mu_i)\mu_i^{-m}((f^n)^*L^{i-1}\cdot L_m\cdot L^{d-i})\\
	&=(D'/\mu_i)h(n,m).
\end{flalign}
Taking $D:=\max\{D'^{i-1}, D'/\mu_i\}$ we get condition (ii) in Lemma \ref{lemexppoly}.
By the computation of mixed degrees in \cite[Corollary 3.4]{Xie2024}, $h$ satisfies (iii) for $\la=\la_{i-1}.$

\medskip

There is $(n,m_0)\in \sN$ 
such that $m_0\geq 2\times 3^bBC^2$ and 
$(n, 2m_0-1)\subseteq \sN.$
By (\ref{equhineqwilcon}) and Lemma \ref{lemexppoly},
we have 
\begin{equation}\label{equationpolympo}\sum_{j=m_0}^{2m_0-1}j^b\leq BC^2(2m_0)^b+m_0^b)\leq BC^22^{b+1}m_0^b
\end{equation}
On the other hand, we have 
$$\sum_{j=m_0}^{2m_0-1}j^b\geq m_0^{b+1}\geq 2\times 3^bBC^2m_0^b.$$
This contradicts (\ref{equationpolympo}). We conclude the proof.
\endproof

\proof[Proof of Lemma \ref{lemexppoly}]
For $p:=(u,v,s,t)\in I:=\C^*\times \C^*\times \Z_{\geq 0}\times \Z_{\geq 0}$, write 
$$H_p:=u^nv^mn^sm^t.$$
For $p:=(u,v,s,t)\in I$,  let $\overline{p}:=(\overline{u},\overline{v},s,t)\in I$ be its complex conjugation.

There is a finite set $S$ of $I$ such that $$h(n,m)=\sum_{p\in S} c_pH_p(n,m)$$ for every $n,m\in \Z_{\geq 0}.$
We may assume further that $c_p\neq 0$ for every $p\in S.$
By assumption (i), after replacing $h$ by $(h+\overline{h})/2$,  we may assume that 
$S$ is invariant under complex conjugacy and we have $c_{\overline{p}}=\overline{c_p}.$
We may assume that $h$ is not identically $0$, then $S\neq \emptyset.$

%

Consider $J:=\R_{>0}\times \R_{>0}\times \Z_{\geq 0}\times \Z_{\geq 0}$  with the lexicographical order $\geq$. 
This is a total order. 
For $p=(u,v,s,t)\in I$, define $|p|:=(|u|, |v|,s,t)\in J.$
Pick $\epsilon_0\in(0,1)$ sufficiently close to $0$ and set $$\sN:=\{(n,m)\in \Z_{\geq 0}|\,\,  (\log n)^2 \leq m\leq n^{\epsilon_0}\}$$ and $$\sN':=\{(n,m)\in \Z_{\geq 0}|\,\,  2^{-1}(\log n)^2 \leq m\leq 2n^{\epsilon_0}\}.$$
It is clear that for $p,p'\in S$,
$|p|> |p'|$
if and only if 
\begin{equation}\label{equesrestterm}\lim_{m\to \infty}\sup\{\frac{|H_{p'}(n,m)|}{|H_p(n,m)|}|\,\,  (n,m)\in \sN'\}\to 0.
\end{equation}
There is a unique maximal element $(\beta,\gamma,a,b)\in \{|p||\,\, p\in S\}.$ 
Set $S^+:=\{p\in S|\,\, |p|=(\beta,\gamma,a,b)\}$ and $S^-=S\setminus S^+.$ Then $S^+, S^-$ are invariant under complex conjugacy.
Set $h^+:=\sum_{p\in S^+}c_pH_p$ and $h^-:=\sum_{p\in S^-}c_pH_p.$ We have $h=h^++h^-.$
By (\ref{equesrestterm}), there is a function $e:\Z_{\geq 0}\to \R_{>0}$ with $\lim\limits_{m\to \infty} e(m)=0$, such that 
\begin{equation}\label{equhmesti}
	|h^-(n,m)|\leq e(m)\beta^n\gamma^mn^am^b
\end{equation}
for every $(n,m)\in \sN.$

\medskip

For every $p\in S^+$, write $$H_p(n,m)=e^{i\theta_pn}e^{i\phi_pm}\beta^n\gamma^mn^am^b$$ where $\theta_p,\phi_p\in \R/2\pi \Z.$
The map $$p\in S^+\mapsto (\theta_p,\phi_p)\in \R/2\pi \Z\times \R/2\pi \Z$$ is injective.
Set $$C(n,m):=\sum_{p\in S^+}c_pe^{i\theta_pn}e^{i\phi_pm}.$$
We have  $h^+(n,m)=C(n,m)\beta^n\gamma^mn^am^b.$ Using Vandermonde determinant, easy to show that $C(n,m)$ is not constantly equal to $0.$

\medskip

Define $\Theta:\Z_{\geq 0}\to (\R/2\pi \Z)^{S^+}$ to be the map $$\Theta: n\mapsto (n\theta_p)_{p\in S^+}.$$
Define $\Phi:\Z_{\geq 0}\to (\R/2\pi \Z)^{S^+}$ to be the map $$\Phi: m\mapsto (m\phi_p)_{p\in S^+}.$$
Define $R:=\Theta\times \Phi: \Z_{\geq 0}\times \Z_{\geq 0}\to (\R/2\pi \Z)^{S^+}\times (\R/2\pi \Z)^{S^+}$ to be the map
$$R: (n,m)\mapsto ((n\theta_p)_{p\in S^+}, (m\phi_p)_{p\in S^+}).$$
Let $q: (\R/2\pi \Z)^{S^+}\times (\R/2\pi \Z)^{S^+}\to \R$ be the function 
$$(x,y)\mapsto \sum_{p\in S^+}c_pe^{ix}e^{iy}.$$
It is continuous. 
We have $$C(n,m)=q\circ R(n,m)=q(\Theta(n),\Phi(m)).$$
Let $Z_1,Z_2$ be the closures of $\Theta(\Z_{\geq 0})$ and $\Phi(\Z_{\geq 0})$ respectively. 
It is clear that $Z_1\times Z_2$ is the closure of $R(\Z_{\geq 0}^2).$
By Poincar\'e recurrence theorem,  there are strictly increasing sequences $n_i, m_i\in \Z_{\geq 0}, i\geq 0$ such that 
$$\Theta(n_i)\to 0 \text{ and } \Phi(m_i)\to 0.$$

\medskip

\begin{lem}\label{lemsubgroup}The subsets
	$Z_1$ and $Z_2$ are subgroups of $(\R/2\pi \Z)^{S^+}$.
\end{lem}

\begin{lem}\label{lemuniformequd}For every non-empty open subset $U\subseteq Z_1\times Z_2$, there is $B_U\geq 0$ such that for every $z\in Z_1\times Z_2$,
	there is $(n,m)\in \{0,\dots, B_U\}^2$ such that $$z+R(n,m)\in U.$$
\end{lem}

We claim that $Z_1\times Z_2$ is the limit set of $R(\sN).$
Let $V$ be any non-empty open subset of $Z_1\times Z_2$.
Let $B_V$ as in Lemma \ref{lemuniformequd}.
For every $T\geq 0$, there is $(n_0,m_0)\in \sN$ with $(n_0,m_0)+\{0,\dots, B_V\}^2\subseteq \sN$ and $m\geq T.$
By Lemma \ref{lemuniformequd}, there is $(n',m')\in \{0,\dots, B_V\}^2$ such that $R(n_0+n',m_0+m')\in V.$ This implies that claim.

\medskip

We claim that $q\geq 0$ on $Z_1\times Z_2.$ 
As $Z_1\times Z_2$ is the limit set of $R(\sN)$, for every $z\in Z_1\times Z_2$, 
there is a sequence $(n'_i,m'_i)\in \sN$ such that $m_i'\to \infty$ and 
$R(n'_i,m'_i)\to z.$ Then we have 
\begin{equation}\label{equcomputeq}0\leq \lim_{i\to \infty}\frac{h(n'_i,m_i')}{\beta^{n_i'}\gamma^{m_i'}{n_i'}^a{m_i'}^b}=\lim_{i\to \infty}\frac{h_+(n'_i,m_i')}{\beta^{n_i'}\gamma^{m_i'}{n_i'}^a{m_i'}^b}=\lim_{i\to \infty}q(R(n'_i,m'_i))=q(z).
\end{equation}
This implies the claim.

As $q$ is not constantly zero on $Z_1\times Z_2$, there is $A_1>0$ and a non-empty open subset $U$ of $Z_1\times Z_2$ such that $q|_{U}>A_1.$
Let $B_U$ as in Lemma \ref{lemuniformequd}.

\begin{lem}\label{lemlowerboundq}
	For every $z\in Z_1\times Z_2$, $$q(z)>\min\{1,\beta\}^{B_U}\min\{\gamma,\beta\}^{B_U}D^{-2B_U}A_1.$$
\end{lem}

\medskip

Set $$A_2:=\min\{1,\beta\}^{B_U}\min\{\gamma,\beta\}^{B_U}D^{-2B_U}A_1$$
and $$A_3:=\max q(Z_1\times Z_2).$$
Pick $T_0\geq 0$ such that $|e(m)|<0.1A_2$ for every $m\geq T_0.$
Then for every $(n,m)\in \sN'$ with $m\geq T_0$, we have 
$$h(n,m)=h^+(n,m)+h^-(n,m)\leq (A_3+0.1A_2)\beta^n\gamma^mn^am^b$$
and 
$$h(n,m)=h^+(n,m)+h^-(n,m)\geq 0.9A_2\beta^n\gamma^mn^am^b.$$

As $h(n,m)>0$ for every $n,m\geq 0$, there is  $C>0$ such that 
\begin{equation}\label{equbetagamnm}C^{-1}\leq \frac{h(n,m)}{\beta^n\gamma^mn^am^b}\leq C
\end{equation}
for $(n,m)\in \sN.$

\medskip

We now assume that (iii) holds and prove the last statement.

Now assume that (iii) holds:
there is $\la>0$ such that for every $\delta\in (0,1)$, there is $D_{\delta}>1$ such that for every $n\geq m\geq 0$, we have 
\begin{equation}\label{equbetanogm}D_{\delta}^{-1}\delta^n\leq \frac{h(n,m)}{\la^n}\leq D_{\delta}\delta^{-n}.
\end{equation}
For $n\gg 0$, $(n,\lfloor (\log n)^2\rfloor+1)\in \sN.$ By (\ref{equbetanogm}),  we have 
$$\lim_{n\to \infty}h(n,\lfloor (\log n)^2\rfloor+1\rfloor)^{1/n}=\la.$$
On the other hand, by (\ref{equbetagamnm}),
$$\lim_{n\to \infty}h(n,\lfloor (\log n)^2\rfloor+1\rfloor)^{1/n}=\beta.$$
Then we get $\beta=\la.$

\medskip

We only need to show that $\gamma=1.$
Let $$S_0:=\{(u,v,s,t)\in S|\,\, |u|=\beta, |v|=\gamma\}.$$
Then $S_0$ is invariant under complex conjugacy and $S^+\subseteq S_0.$
Set $C_0:=\max\{s+t| (u,v,s,t)\in S\}.$
Set $$S':=\{(u,v,s,t)\in S_0|\,\, s+t=C_0\}.$$
Then $S'$ is non-empty and  invariant under complex conjugacy.
Set $S'':=S\setminus S'.$
For $\epsilon_1\in (0,1)$, set $$\sM(\epsilon_1):=\{(n,m)|\,\, 2^{-1}\epsilon_1n\leq m\leq \epsilon_1 n\}.$$
It is clear that there is $\epsilon'\in (0,1)$ such that for every $\epsilon_0\in (0,\epsilon')$ and
$p\in S'$, we have 
\begin{equation}\label{equesresttermsp}\lim_{m\to \infty}\sup\{\frac{|H_{p}(n,m)|}{\beta^n\gamma^mn^{C_0}}|\,\,  (n,m)\in \sM(\epsilon_1)\}\to 0.
\end{equation}

Set $h':=\sum_{p\in S'}c_pH_p$ and $h'':=\sum_{p\in S''}c_pH_p.$ We have $h=h'+h''.$
By (\ref{equesresttermsp}), there is a function $e':\Z_{\geq 0}\to \R_{>0}$ with $\lim\limits_{m\to \infty} e'(m)=0$, such that 
\begin{equation}\label{equhmestihp}
	|h'(n,m)|\leq e'(m)\beta^n\gamma^mn^{C_0}
\end{equation}
for every $(n,m)\in \sN.$
For every $p\in S'$, write $$H_p(n,m)=e^{i\theta_pn}e^{i\phi_pm}\beta^n\gamma^mn^{C_0}(\frac{m}{n})^{b_p}$$ where $\theta_p,\phi_p\in \R/2\pi \Z.$
Set $$D(n,m):=\sum_{p\in S'}c_pe^{i\theta_pn}e^{i\phi_pm}(\frac{m}{n})^{b_p}.$$
We have  $h'(n,m)=D(n,m)\beta^n\gamma^mn^{C_0}.$ 
Set $b_-:=\min\{t|\,\, (u,v,s,t)\in S' \}.$ Set $S'_+:=\{(u,v,s,t)\in S'|\,\, t=b_-\}.$
So $S'_+$ is non-empty and invariant under complex conjugacy.
Set $S'_-:=S'\setminus S'_+.$
Set $$C_+(n,m):=\sum_{p\in S'_+}c_pe^{i\theta_pn}e^{i\phi_pm}$$
and $$D_+(n,m):=C_+(n,m)(\frac{m}{n})^{b_-}.$$
Set $$D_-(n,m):=\sum_{p\in S'_-}c_pe^{i\theta_pn}e^{i\phi_pm}(\frac{m}{n})^{b_p}$$
Set $$A_-:=\sum_{p\in S'_-}|c_p|.$$
Then for $(n,m)\in \sM(\epsilon_1), $
\begin{equation}\label{equupdm}
	|D_-(n,m)|\leq A_-\epsilon_1^{b_-+1}.
\end{equation}

Define $\Theta':\Z_{\geq 0}\to (\R/2\pi \Z)^{S'_+}$ to be the map $$\Theta': n\mapsto (n\theta_p)_{p\in S'_+}.$$
Define $\Phi':\Z_{\geq 0}\to (\R/2\pi \Z)^{S'_+}$ to be the map $$\Phi: m\mapsto (m\phi_p)_{p\in S'_+}.$$

Define $R':=\Theta\times \Phi: \Z_{\geq 0}\times \Z_{\geq 0}\to (\R/2\pi \Z)^{S'_+}\times (\R/2\pi \Z)^{S'_+}$ to be the map
$$R': (n,m)\mapsto ((n\theta_p)_{p\in S'_+}, (m\phi_p)_{p\in S'_+}).$$
Let $q': (\R/2\pi \Z)^{S'_+}\times (\R/2\pi \Z)^{S'_+}\to \R$ be the function 
$$(x,y)\mapsto \sum_{p\in S'_+}c_pe^{ix}e^{iy}.$$
It is continuous. 
We have $$C_+(n,m)=q'\circ R'(n,m)=q'(\Theta'(n),\Phi'(m)).$$
Note that the map $$p\in S'_+\mapsto (\theta_p,\phi_p)\in \R/2\pi \Z\times \R/2\pi \Z$$ is injective.
Using Vandermonde determinant, easy to show that $C'(n,m)$ is not constantly equal to $0.$
In particular $q'$ is not constantly equal to $0$ on $Z_1'\times Z_2'.$

Let $Z_1',Z_2'$ be the closures of $\Theta'(\Z_{\geq 0})$ and $\Phi(\Z_{\geq 0})$ respectively. 
Same proof as Lemma \ref{lemsubgroup} shows that $Z_1',Z_2'$ are closed subgroups in $(\R/2\pi \Z)^{S'_+}.$
There is a non-empty open subset $U'\subseteq Z_1'\times Z_2'$  and a constant $A_3>0$ such that for every $z\in U'$, $|q'(z)|>A_3.$
For every $(n,m)\in \sM(\epsilon_1)$, if $R(n,m)\in U'$, then 
\begin{equation}\label{equbounddpinup}
	|D_+(n,m)|\geq 2^{-b_-}A_3\epsilon_1^{b_-}
\end{equation}
Let $A_4:=\max\{|q(z)||\,\, z\in Z_1'\times Z_2'\}.$
Then 
\begin{equation}\label{equbounddpinupup}
	|D_+(n,m)|\leq A_4\epsilon_1^{b_-}
\end{equation}

Fix $\epsilon_1\in (0, \epsilon')$ with $\epsilon_1<\frac{0.1A_3}{A_12^{b_-}}.$
Then for every $(n,m)\in \sM(\epsilon_1)$, if $R(n,m)\in U'$, we have 
\begin{equation}\label{bounddinup}|D(n,m)|\geq 0.9\times 2^{-b_-}A_3\epsilon^{b_-} \text{ and } |D(n,m)|\leq 2A_3\epsilon^{b_-}.
\end{equation}
Same proof as Lemma \ref{lemuniformequd} shows that there is $B'\geq 0$ such that for every $z\in Z_1'\times Z_2'$,
there is $(n,m)\in \{0,\dots, B'\}^2$ such that $$z+R'(n,m)\in U'.$$
For every $n\geq 0$, there is $(N(n),M(n))\in (\lfloor 0.75\epsilon n\rfloor, n)+\{0,\dots, B'\}^2$ such that 
$$R'(N(n),M(n))\in U'.$$
By (\ref{bounddinup}) and (\ref{equhmestihp}), there is $C_1>1$ such that for every $(n,m)\in \sM(\epsilon_1)$, we have
\begin{equation}\label{equmeesti}
	C_1^{-1}<\frac{h(N(n),M(n))}{\beta^{N(n)}\gamma^{M(n)}N(n)^{C_0}}<C_1.
\end{equation}
Then we have 
$$\lim_{n\to \infty}h(N(n),M(n))^{1/n}=\lim_{n\to \infty}(\beta^{N(n)}\gamma^{M(n)}N(n)^{C_0})^{1/n}=\beta\times \gamma^{\epsilon_0}.$$
By (iii) and the fact $\beta=\la$, we get 
$$\lim_{n\to \infty}h(N(n),M(n))^{1/n}=\la.$$
Then we get $\gamma=1$, which concludes the proof.
\endproof

\proof[Proof of Lemma \ref{lemsubgroup}]
As $\Theta$ and $\Phi$ are group homeomorphisms, both $Z_1$ and $Z_2$ contain $0$ and  are closed under addition. 
For every $z\in Z_1$, there is an increasing sequence $l_i\in \Z_{\geq 0}$ such that $$\Theta(l_i)\to z.$$
Note that $l_i$ may not be strictly increasing.
After taking subsequence, we may assume that $n_i\geq l_i.$ Then we have $$\Theta(n_i-l_i)\to -z.$$
This implies that $Z_1$ is a subgroup. Similarly, $Z_2$ is a subgroup. 
\endproof

\proof[Proof of Lemma \ref{lemuniformequd}]
Pick $w\in U$. There is an open neighborhood $W$ of $0$ such that $W=-W$ and $W+W\subseteq U-w.$ Set $V:=w+W\subseteq U.$ 
We have $W+V\subseteq U.$
As $Z_1\times Z_2$ is a group, for every $z\in Z_1\times Z_2$, $-z\in Z_1\times Z_2$.
So there is $(n_z,m_z)\in \Z_{\geq 0}^2$ such that $R(n_z,m_z)\in V-z.$
As $Z_1\times Z_2$ is compact, there is a finite subset $F$ of $Z_1\times Z_2$ such that 
$F+V=Z_1\times Z_2.$ For every $z'\in Z_1\times Z_2$, there is $z''\in F$ such that $z'\in z''+W.$
Then we have 
$$z'+R(n_{z''},m_{z''})=z''+R(n_{z''},m_{z''})+(z'-z'')\subseteq V+W\subseteq U.$$
Setting $B_U:=\max\{\max\{n_z, m_z\}| z\in F\},$ we conclude the proof.
\endproof	

\proof[Proof of Lemma \ref{lemlowerboundq}]
Pick a sequence $(n'_i,m'_i)\in \sN$ such that $m_i'\to \infty$, and  $R(n_i',m_i')\to z.$ 
There is $(b_1,b_2)\in \{0,\dots, B_U\}^2$ such that $$R(b_1,b_2)+z\in U.$$
Note that for $m$ sufficiently large, we have $\sN+\{0,\dots, B_U\}^2\subseteq \sN'.$

By (\ref{equcomputeq}), we have 
$$\lim_{i\to \infty}\frac{h(n'_i,m_i')}{\beta^{n_i'}\gamma^{m_i'}{n_i'}^a{m_i'}^b}=q(z).$$
Similarly, we have 
$$\beta^{-b_1}\gamma^{-b_2}\lim_{i\to \infty}\frac{h(n'_i+b_1,m_i'+b_2)}{\beta^{n_i'}\gamma^{m_i'}{n_i'}^a{m_i'}^b}=\lim_{i\to \infty}\frac{h(n'_i+b_1,m_i'+b_2)}{\beta^{n_i'}\gamma^{m_i'}{(n_i'+b_1)}^a{(m_i'+b_2)}^b}$$
$$=q(R(b_1,b_2)+z)>A_1.$$
Our condition (ii) implies that 
$$h(n'_i+b_1,m_i'+b_2)\leq D^{b_1+b_2}h(n_i',m_i')\leq D^{2B_U}h(n_i',m_i').$$
The we get 
$$A_1<\beta^{-b_1}\gamma^{-b_2}D^{2B_U}\lim_{i\to \infty}\frac{h(n'_i,m_i')}{\beta^{n_i'}\gamma^{m_i'}{n_i'}^a{m_i'}^b}=\beta^{-b_1}\gamma^{-b_2}D^{2B_U}q(z).$$
So we have $$q(z)>\beta^{b_1}\gamma^{b_2}D^{-2B_U}A_1\geq \min\{1,\beta\}^{B_U}\min\{\gamma,\beta\}^{B_U}D^{-2B_U}A_1,$$
which concludes the proof.
\endproof

\begin{cor}\label{corfactoramplifiedquasi}Let $Y$ be a projective variety over $\bk$ and $g:Y\to Y$ be an endomorphism. Let $\pi: X\to Y$ be a surjective morphism such that $\pi\circ f=g\circ \pi.$ If $f$ is $\alpha$-quasi-amplified for some $\alpha\in \R_{>0}$, then $g$ is $\alpha$-quasi-amplified.
\end{cor}
\proof
The product formula for relative dynamical degrees (c.f. \cite{Dinh2011}, \cite{Dang2020} and \cite[Theorem 1.3]{Truong2020}) shows that $$\{\mu_i(V,g)|\,\, i=1,\dots, d_V\}\subseteq \{\mu_i(W,f)|\,\, i=1,\dots, d_W\}.$$
We conclude the proof by Theorem \ref{thmchendamplify}.
\endproof

\section{Measures for the constructible topology}
Let $X$ be a reduced projective scheme over $\bk$ of dimension $d$. For every $x\in X$, denote by $Z_x:=\overline{\{x\}}.$

\subsection{Constructible topology}
Denote by $|X|$ the underling set of $X$ with the constructible topology; i.e. the topology on a  $X$ generated by the constructible subsets (see~\cite[Section~(1.9) and in particular (1.9.13)]{EGA-IV-I}).
In particular every constructible subset is open and closed.
This topology is finer than the Zariski topology on $X.$ Moreover $|X|$ is (Hausdorff) compact. 

\medskip

Let $\sC(|X|)$ be the space of continuous functions on $|X|$ endowed with the norm $\|h\|:=\max\{f(x)|\,\, x\in |X|\}.$
Let $\sA(|X|)$ be the set of constructible subsets of $X$. For every $U\in \sA(|X|),$ denote by $1_U$ the characteristic function of $U$ i.e. 
$1_U(x)=1$ if $x\in U$ and $1_U(x)=0$ if $x\not\in U.$  Such functions are continuous, as constructible subsets are both open and closed in $|X|$. For two distinct points $x,y\in |X|$, there is a constructible subset $U$ of $X$ such that $x\in U$ and $y\not\in U.$
So we have $1_U(x)\neq 1_U(y).$
By Stone-Weierstrass theorem, we get the following result.
\begin{lem}\label{lemdensesetcon}The $\R$-algebra generated by $\{1_U|\,\, U\in \sA(|X|)\}$ is dense in $\sC(|X|)$.
\end{lem}

\subsection{Measures}
Let $\sA(|X|)$ be the set of constructible subsets of $X$. It is a field of sets in $|X|$ i.e. it has the following properties:
\begin{points}
	\item $\emptyset\in \sA(|X|), |X|\in \sA(|X|)$;
	\item for $U,V\in \sA(|X|)$, $U\cap V\in \sA$ and $U\setminus V\in \sA(|X|).$
\end{points}
Moreover, the Borel $\sigma$-algebra $\sB(|X|)$ on $|X|$ is generated by $\sA(|X|)$.

\begin{pro}\label{proxconradon}Every finite measure Borel measure on $|X|$ is a Radon measure i.e. $|X|$ is a Radon space.
\end{pro}
\proof
Let $\mu$ be a finite Borel measure on $|X|$. By Riesz-Markov-Kakutani representation theorem (c.f. \cite[Theorem 2.14]{Rudin1987}), there is a Radon measure $\mu'$ such that 
for every continuous function $h$ on $|X|$, $\int h\mu=\int h\mu'.$ 
So for every constructible subset $V$ of $X$, we have $$\mu'(V)=\int 1_V\mu'=\int 1_V\mu=\mu(V).$$
As the restriction of $\mu$ and $\mu'$ on $\sA$ are the same and both of them are finite, the uniqueness part of the 
Carath\'eodory's extension theorem, implies that $\mu=\mu'$. So $\mu$ is a Radon measure, which concludes the proof.
\endproof

Denote by $\sM(|X|)$ the space of finite Borel measures on $X$ endowed with the weak-$\ast$ topology.
Combining \cite[Theorem 1.12]{Xie2023}  with Proposition \ref{proxconradon}, we get the following result.
\begin{thm}\label{thmRadon}Every $\mu\in \sM(|X|)$ takes form 
	$$\mu=\sum_{i\geq 0}a_i\delta_{x_i}$$
	where $\delta_{x_i}$ is the Dirac measure at $x_i\in X$, $a_i\geq 0$ with $\sum_{i\geq 0}a_i<+\infty$.
\end{thm}
In \cite[Theorem 1.12]{Xie2023}, $X$ is assumed to be irreducible. However the proof still work without the irreducibility.
Let $\sM^1(|X|)$ be the space of probability measures on $|X|.$ 
Since $|X|$ is compact, $\sM^1(|X|)$ is compact.  
Combining \cite[Theorem 1.12]{Xie2023}  with Proposition \ref{proxconradon}, we get the following result.
\begin{cor}\label{corsmxsc}
	The space $\sM^1(|X|)$ is sequentially compact.
\end{cor}
In \cite[Theorem 1.12]{Xie2023}, $X$ is assumed to be irreducible. However the proof still works without the irreducibility.

\subsection{Vector valued measures}
Let $V$ be a finitely dimensional $\R$-space. Let $\sP$ be a convex closed salient cone in $V$ with $\sP^{\circ}\neq\emptyset.$ Let $\sP^{\vee}$ be its dual cone in $V^{\vee}$. Then $\sP^{\vee}$ is a convex closed salient cone in $V^{\vee}$ with $(\sP^{\vee})^{\circ}\neq\emptyset.$ Let $\|\cdot\|$ be any norm on $V.$

\begin{rem}\label{remconenorm}
	For every $\omega\in (\sP^{\vee})^{\circ}$, we may define a norm $\|\cdot\|_{\omega}$ on $V$ as follows:
	for every $v\in V$, $\|v\|_{\omega}:=\inf\{((v_1+v_2)\cdot \omega)|\,\,v_1,v_2\in \sP, v_1-v_2=v\}.$
	Easy to check that $\|\cdot\|_{\omega}$ is a norm. Moreover, for every $v\in \sP$, we have $\|v\|_{\omega}=(v\cdot \omega).$
\end{rem}

\medskip

We define a $\sP$-valued Borel measure on $|X|$ to be a map $\gamma: \sB(|X|)\to \sP$ such that 
\begin{points}
	\item $\gamma(\emptyset)=0;$
	\item for a countable collection of disjoint Borel subsetes $E_n, n\geq 0$, we have 
	$$\gamma(\sqcup_{n\geq 0}E_n)=\sum_{n\geq 0}\gamma(E_n).$$
\end{points}

Denote by $\sM(|X|,\sP)$ the set of $\sP$-valued Borel measure on $|X|$. 
There is a map $$(\beta,\gamma)\in\sP^{\vee}\times \sM(|X|,\sP)\to (\beta\cdot \gamma)\in \sM(|X|)$$ 
such that for every $E\in \sB(|X|)$, $(\beta\cdot \gamma)(E)=(\beta\cdot \gamma(E)).$
Easy to check that $(\beta\cdot \gamma)$ is a Borel measure on $|X|.$
Define the weak-* topology on $\sM(|X|,\sP)$ by the topology generated by the sets of form 
$$\{\gamma\in \sM(|X|,\sP)|\,\, \int h (\beta\cdot \gamma)<c\}$$ where $\beta\in \sP^{\vee}, h\in C(|X|)$ and $c\in \R.$
\medskip

On the other hand, once we have a continuous map  $$\phi:\sP^{\vee}\to \sM(|X|)$$ such that for $r_1,r_2\geq 0, v_1,v_2\in \sP^{\vee}$, we have
$$\phi(r_1v_1+r_2v_2)=r_1\phi(v_1)+r_2\phi(v_2),$$
then $\phi$ defined a $\sP$-valued Borel measure $\nu_{\phi}\in \sM(|X|,\sP)$ such that for every $E\in \sB(|X|)$,
$\nu(E)$ is the unique element in $\sP$ such that for every $v\in \sP^{\vee}$, $$(v,\nu(B))=\phi(v)(E).$$

\medskip

Now we give a concrete description of $\sM(|X|,\sP).$
Let $\sM^{\pm}(|X|)=\sM(|X|)\otimes_{\R_{>0}}\R$ be the space of signed measures.
By Theorem\ref{thmRadon}, every $\mu\in \sM^{\pm}(|X|)$ takes form 
$$\mu=\sum_{i\geq 0}a_i\delta_{x_i}$$
where $\delta_{x_i}$ are the Dirac measures at distinct points $x_i\in X$ and $a_i\in \R$ with $\sum_{i\geq 0}|a_i|<+\infty$.
Consider the $\R$-space $\sM(|X|)\otimes_{\R_{>0}}V:=\sM^{\pm}\otimes_{\R}V.$
Assume that $\dim V=s$ and fix $\beta_1^{\vee},\dots, \beta_s^{\vee}\in (\sP^{\vee})^{\circ}$ which forms a base of $V^{\vee}$.
Let $\beta_1,\dots, \beta_s$ be its dual basis in $V$. Then $\sM(|X|)\otimes_{\R_{>0}}V=\oplus_{i=1}^s\sM^{\pm}(|X|)\beta_i.$
So every element in $\sM(|X|)\otimes_{\R_{>0}}V$ takes form 
$$\mu=\sum_{i\geq 0}\alpha_i\delta_{x_i}$$
where $\delta_{x_i}$ are the Dirac measures at distinct points $x_i\in X$ and $\alpha_i\in V$ with $\sum_{i\geq 0}\|\alpha_i\|<+\infty$.

Consider the map $$\Theta:\gamma\in \sM(|X|,\sP)\mapsto (\beta_i^{\vee}\cdot \gamma)\beta_i \in \sM(|X|)\otimes_{\R_{>0}}V.$$
This map induces a homeomorphism from  $\sM(|X|,\sP)$ to its image. The pairing $\sP^{\vee}\times \sM(|X|,\sP)\to \sM(|X|)$ extends to 
the pairing $$\sM(|X|)\otimes_{\R_{>0}}V\times V^{\vee}\to \sM^{\pm}(|X|)$$ defined in the obvious way:
$$((\oplus_{i=1}^s\mu_i\beta_i) \cdot (\sum_{i=1}^sa_i\beta_i^{\vee})):=\sum_{i=1}^sa_i\mu_i.$$ 
Then the image of $\Theta$ is the intersection $$\Theta(\sM(|X|,\sP))=\cap_{\beta\in \sP^{\vee}}\{\alpha\in \sM(|X|)\otimes_{\R_{>0}}V|\,\,(\alpha\cdot \beta)\in \sM(|X|)\}.$$

\medskip

For every Borel set $E\in \sB(|X|)$, the characteristic function $1_E$ is measurable. 
In particular, for every $E\in \sB(|X|)$, and $\gamma\in \sM(|X|,\sP)$, we get $1_E\gamma \in \sM(|X|,\sP).$ 
If further $E\in \sA(|X|)$, $1_E$ is continuous. So the endomorphism $$\gamma\in\sM(|X|,\sP)\to 1_E\gamma\in \sM(|X|,\sP)$$ is continuous. 

\medskip

For every $x\in |X|$ and $\alpha\in \sP$, define $\alpha\delta_x\in \sM(|X|,\sP)$ as follows: for every $E\in \sB(|X|)$, $\alpha\delta_x(E)=\alpha$ if $x\in V$ and $\alpha\delta_x(E)=0$ if $x\not\in E.$ More generally, for an at most countable subset $S$ of $|X|$, with a sequence of vectors $\alpha_x\in \sP, x\in S$ with $\sum_{x\in S}\|\alpha_x\|<+\infty,$ we may define $\sum_{x\in S}\alpha_x\delta_x\in  \sM(|X|,\sP)$ in the obvious way: for every $V\in \sB(|X|)$, $$(\sum_{x\in S}\alpha_x\delta_x)(E):=\sum_{x\in S}\alpha_x\delta_{x}(E).$$
Easy to check that $$\Theta(\sum_{x\in S}\alpha_x\delta_x)= \sum_{x\in S}\alpha_x\delta_x\in \sM(|X|)\otimes_{\R_{>0}}E.$$

On the other hand, for every $\gamma\in \sM(|X|,\sP)$, we may write 
$\Phi(\gamma)=\sum_{x\in S}\alpha_x\delta_x$ for an at most countable subset $S$ of $|X|$, $\alpha_x\in V$ with $\sum_{x\in S}\|\alpha_x\|<+\infty.$ Then for every $\beta\in \sP^{\vee}$, $$(\beta\cdot \gamma)=\sum_{x\in S}(\beta\cdot \alpha_x)\delta_x\in \sM(|X|).$$ So $(\beta\cdot \alpha_x)\geq 0$ for every $x\in S.$ So we have $\alpha_x\in \sP.$
We then get the following generalization of Theorem \ref{thmRadon}.
\begin{pro}\label{provectmeasure}Every $\gamma\in \sM(|X|,\sP)$ takes form 
	$$\gamma=\sum_{x\in S}\alpha_x\delta_x$$ for an at most countable subset $S$ of $|X|$, $\alpha_x\in \sP$ with $\sum_{x\in S}\|\alpha_x\|<+\infty.$
\end{pro}

For $\gamma\in \sM(|X|,\sP)$, the support $\Supp\,\gamma$ of $\gamma$ is defined to be the support of the measure $(\beta\cdot \gamma)$ where $\beta\in (\sP^{\vee})^{\circ}.$ It does not depend on the choice of $\beta$. Indeed $|X|\setminus \Supp\,\gamma$ is the maximal open subset $U$ of $|X|$ satisfying $\gamma(U)=0.$ 

\medskip

Let $\sM^1(|X|,\sP)$ be the subset elements $\gamma$ of $\sM(|X|,\sP)$ with $\|\gamma(X)\|\leq 1.$ By Proposition \ref{provectmeasure} and Corollary \ref{corsmxsc}, $\sM^1(|X|,\sP)$ is compact and sequentially compact.

\section{Generated cycles}\label{sectiongencycle}
Let $X$ be a reduced projective scheme over $\bk$ of dimension $d$. 
For every $i=0,1,\dots, d$, denote by $X_i$ the set of (scheme-theoretic) points $x\in X_i$ with $\dim Z_x=i$ where $Z_x:=\overline{\{x\}}.$
Denote by $Z_i(X)_{\R}$ the space of $i$-cycles with $\R$-coefficients. 
Every $Z\in Z_i(X)_{\R}$ can be uniquely written as $$Z=\sum_{x\in X_i} m(Z,x)Z_x,$$ where $m(Z,x)=0$ for all but finitely many $x\in X_i.$
Let $c(Z)$ be the set of $x$ with $m(Z,x)\neq 0.$
Denote by $\Eff_i(X)$ the subset of effective $i$-cycles in  $Z_i(X)_{\R}$.
For every subset $U$ of $X$. Denote by $Z_i(U,X)_\R$ the set of $i$-cycles $Z$ with $c(Z)\subseteq U$; and $\Eff_i(U,X)$  the subset of effective $i$-cycles in  $Z_i(U,X)_\R$.
Denote by $\Psef_i(X)$ the cone of pseudo-effective classes in $N_i(X)_{\R}.$

\medskip

Let $U$ be a locally closed subset of $X$ i.e. $U$ is open in its Zariski closure $\overline{U}^{\zar}$.
Let $\Psef^i(U,X)$ be the closure of the convex cone in $\overline{U}^{\zar}$ generated by effective $i$-cycles of form $$\sum_{j=1}^ma_jZ_{x_j}$$ where $a_j\geq 0$, $x_j\in U.$ For $U_1\subseteq U_2$ with the same Zariski closure, we have $\Psef^i(U_1,X)\subseteq \Psef^i(U_2,X).$

\medskip

\subsection{Generated cycles}
Let $\sZ(X)$ be the set of all non-empty Zairski closed subsets of $X$.
Moreover, for Zariski closed subsets $V_1\subseteq V_2$, we have natural map $\iota_{V_1\subseteq V_2}: N_i(V_1)_\R\to N_i(V_2)_\R$ induced by the inclusion $V_1\hookrightarrow V_2$. {\bf For simplifying the notations, we often omit the morphism $\iota_{V_1\subseteq V_2}$.}
We note that if $\dim V<i$,  then $N_i(V)_{\R}=0.$
For every $V\in \sZ(X)_{\R}$, set $R_V: Z_i(X)_{\R}\to N_i(V)_{\R}$ sending $Z$ to $$R_V(Z):=\sum_{x\in X_i\cap V} m(Z,x)[Z_x].$$
Define $$\Phi:=\prod_{V\in \sZ(X)}R_V: Z_i(X)_\R\to \prod_{V\in \sZ(X)}N_i(V)_\R.$$
It is clear the $\Phi$ is injective.  
We endow $\prod_{V\in \sZ(X)}N_i(V)_\R$ the product topology. 
Define $\sG_i(X)$ to be the closure of $\Phi(Z_i(X)_\R)$ with the induced topology. We call elements in $\sG_i(X)$ the \emph{generated $i$-cycles}.
We now identify $Z_i(X)_\R$ with its image in $\sG_i(X).$
We would like to think generated $i$-cycles as an analogy of the notion of closed current of bidimension $(i,i)$ in complex geometry. 
For every $V\in \sZ(X)$, the morphism $R_V$ is just the restriction of the projection $\prod_{W\in \sZ(X)}N_i(W)_\R\to N_i(V)_{\R}$ to $Z_i(X)_R$. So $R_V$ extends to a continuous morphism $$R_V: \sG_i(X)\to N_i(V)_{\R}.$$
The continuity of $R_V, V\in \sZ(X)$ implies the following cut-and-paste relations: for every $V_1, V_2\in \sZ$, we have 
\begin{equation}\label{equationcutpaste}R_{V_1}+R_{V_2}=R_{V_1\cup V_2}+R_{V_1\cap V_2}
\end{equation}
in $N_i(V_1\cup V_2)_{\R}.$

\begin{lem}\label{lemdescpgs}The following holds:
\begin{points}
\item $\sG_i(X)=\{\alpha\in \prod_{V\in \sZ(X)}N_i(V)_\R|\,\, (\ref{equationcutpaste}) \text{ holds at } \alpha\}$.
\item the restriction of the projection $\Psi: \prod_{V\in \sZ(X)}N_i(V)_{\R}\to \prod_{x\in X}N_i(Z_x)_{\R}$ induces an isomorphism 
$$\Psi|_{\sG(X)}: \sG(X)\simeq \prod_{x\in X}N_i(Z_x)_{\R}.$$
\end{points}
\end{lem}
\proof
Set $H:= \{\alpha\in \prod_{V\in \sZ(X)}N_i(V)_\R|\,\, (\ref{equationcutpaste}) \text{ holds at } \alpha\}.$ By (\ref{equationcutpaste}),
we have $\sG_i(X)\subseteq H.$

Next we show that $\Psi|_H: H\to \prod_{x\in X}N_i(Z_x)_{\R}$ is an isomorphism. For this, only need to construct an inverse $h$ of $\Psi|_H$.
Consider $\Z_{\geq 0}\times \Z_{>0}$ with lexicographical order $\geq$ which is a well ordered set. 
For every $V\in \sZ(X)$, let $s_V$ be the number of irreducible components of $V$.
Define $d(V):=(\dim V, s_V).$ For $\alpha:=(\alpha_x)_{x\in |X|}\in \prod_{x\in X}N_i(Z_x)_{\R}$, we define $h(\alpha)_V$ by induction on $d(V)$.
If $s_V=1$, let $\eta$ be the generic point of $V$. Define $h(\alpha)_V:=\alpha_{\eta}.$
In particular, $h(\alpha)_V$ defined when $d(V)=(0,1).$
Now assume that $s_V>1$. Let $V_1,\dots, V_s$ by the generic points of the irreducible components of $V.$
Note that $s\geq 2$. Set $W:=V_2\cup\dots\cup V_s.$
We have $$\max\{d(V_1\cap W),d(V_1), d(W)\}< d(V)$$
Define $$h(\alpha_V):=h(\alpha)_{V_1}+h(\alpha)_{W}-h(\alpha)_{V_1\cap W}$$ in $N_i(V)_{\R}.$
Easy to check that $h(\prod_{x\in X}N_i(Z_x)_\R)\subseteq H$, $\Psi\circ h=\id$, $h\circ \Psi|_{H}=\id$ and
$h$ is continuous. This implies that $\Psi|_H: H\to \prod_{x\in X}N_i(Z_x)_{\R}$ is an isomorphism.

\medskip

We now only need to show that $\Psi(Z_i(X)_{\R})$ is dense in $\prod_{x\in X}N_i(Z_x)_\R.$  
We only need to show that for every  $x_1,\dots, x_s\in X$ and $\alpha_l\in N_i(Z_{x_l})_\R$, there is $\alpha\in Z_i(X)_\R$ such that 
$\Psi(\alpha)_{x_l}=\alpha_l.$ For every $l=1,\dots, s$, let $D(l)$ be the set of index $j$ such that $x_j\in Z_{x_l}.$
Set $W_l:=Z_{x_l}\cap (\cup_{j\neq l} Z_{x_j}).$
\begin{lem}\label{lemfreeclass}For every  $l=1,\dots, s$, there is $\beta_l\in Z_i(Z_{x_l})_{\R}$ such that the image of $\beta_l$ in $N(Z_{x_l})_\R$
is $\alpha_{l}-\sum_{j\in D(l)}\alpha_j$ and we have $R_{W_l}(\beta_l)=0$.
\end{lem}
Easy to check that $\Psi(\sum_{l=1}^s\beta_l)_{j}=\alpha_j$ for every $j=1,\dots, s$. This concludes the proof.
\endproof
\proof[Proof of Lemma \ref{lemfreeclass}]
By De Jong's alteration theorem \cite{Jong1996}, there is a smooth projective variety $Y$ with a generically finite and surjective morphism $q: Y\to V_i$.
Pick any  $\gamma\in Z_i(Z_{x_l})_{\R}$ such that the image of $\gamma$ in $N(Z_{x_l})_\R$
is $\alpha_{l}-\sum_{j\in D(l)}\alpha_j.$ Pick $\gamma'\in Z_i(Y)_{\R}$ such that $q_*\gamma'=\gamma.$
By Chow's moving lemma \cite[Lemma 43.24.1]{stacks-project}, there is $\gamma''\in Z_i(Y)_{\R}$ which is linearly equivalent to $\gamma'$ and no irreducible component of $\Supp\, \gamma'$
contained in $q^{-1}(W_l).$ Then $\beta_l:=\gamma''$ satisfies the conditions we need.
\endproof

Let $g: Y\to X$ be a morphisms between reduced projective schemes over $\bk$.
Easy to check that the pushforward $g_*: Z_i(Y)\to Z_i(X)$ extends to a continuous morphism $g_*: \sG_i(Y)\to \sG_i(X).$

\subsubsection{Restriction map}
For every $V\in \sZ(X)$, consider the restriction map $$\pi_V: \prod_{W\in \sZ(X)}N_i(W)_\R\to \prod_{W'\in \sZ(V)}N_i(W')_\R$$
sending $(\alpha_W)_{W\in \sZ(X)}$ to $(\alpha_{W'\cap V})_{W'\in \sZ(X)}.$ It is clear that $\pi_V$ is continuous. For every $Z\in Z_i(X)$, we have 
$$\pi_V\circ \Phi(Z)=\Phi(\sum_{x\in X_i\cap V} m(Z,x)Z_x).$$
So $\pi_V(\sG_i(X))\subseteq \sG_i(V).$
We still denote by $\pi_V$ its restriction on $\sG_i(X).$
Let $i_V: V\hookrightarrow X$ be the inclusion morphism, then we have 
$$\pi_V\circ (i_V)_*=(i_V)_*.$$
In particular, $(i_V)_*$ gives a natural embedding from $\sG_i(V)$ to $\sG_i(X).$ We identify $\sG_i(V)$ with its image under $(i_V)_*$ in $\sG_i(X).$

\medskip

We now extend the above definition to any constructible subset by induction on its dimension.
Let $W$ be a constructible subset. If $\dim W=0$, then $W$ is closed. So $\pi_W$ is defined.
Now assume that $\dim W=r\geq 1$ and the restriction map is defined for every constructible set of dimension $<r.$
We define $$\pi_W:=\pi_{\overline{W}}-\pi_{\overline{W}\setminus W}.$$
The two terms in the right hand side are defined as $\overline{W}$ is closed and $\dim (\overline{W}\setminus W)<r.$
By induction, easy to check that $\pi_W$ is continuous and for every 
$Z\in Z_i(X)$, we have 
\begin{equation}\label{equresconstcyc}\pi_W\circ \Phi(Z)=\Phi(\sum_{x\in X_i\cap W} m(Z,x)Z_x).
\end{equation}

Then for two disjoint constructible sets $W_1,W_2$, we get 
\begin{equation}\label{equationdisjdec}\pi_{W_1\sqcup W_2}=\pi_{W_1}+\pi_{W_2} \text{ and } \pi_{W_1}\circ \pi_{W_2}=0.
\end{equation}
Then for constructible sets $W, W'$ with $W\subseteq W'$, we have $$\pi_W=\pi_W\circ \pi_W'.$$

Define $\sG_i(W,X):=\pi_W(\sG_i(X)).$ By (\ref{equationdisjdec}), $\sG_i(W,X)$ is the closure of $Z_i(W,X)_{\R}$ in $\sG_i(X).$

\subsubsection{Support of generated cycles}
Let $Z\in \sZ_i(X)$. For every $x\in |X|$, we say that $x$ is \emph{{\bf not} in the support of $Z$} if there is a constructible subset $W$ containing $x$ such that $\pi_W(Z)=0.$
Otherwise, we say that $x$ is in the support of $Z$.
Define $\Supp\, Z$ to be the set of all $x\in |X|$ which is in the support of $Z$. As constructible subsets are open in $|X|$, $\Supp\, Z$ is closed in $|X|.$
By compactness of constructible subsets, for every constructible subset $W\subseteq |X|\setminus \Supp\, Z$, $\pi_W(Z)=0.$

\medskip

For every closed subset $Y$ of $|X|$, we say that $Z$ is supported on $Y$ if $\Supp\, Z\subseteq Y.$
If $Y$ is a constructible subset, then $\sR_i(Y,X)=\pi_Y(\sR_i(X))$ is exactly the 
 subspace of $Z\in \sR_i(X)$ supported on $Y$.  

Let $X=\sqcup_{i=1}^sY_i$ be a finite partition of $X$ by constructible subsets $Y_i.$ Then we get a direct composition 
$$\sG_i(X)=\oplus_{i=1}^s \sG_i(Y, X).$$

\subsubsection{Intersection numbers}
For every $\alpha\in \sG_i(X)$, $Z\in \sZ(X)$ and $\beta\in N^i(Z)_{\R}$, define $$(\alpha\cdot \beta):=(R_Z(\alpha)\cdot \beta).$$
It is clear that for every $Z\in \sZ(X)$, the map $$(\alpha,\beta)\in \sG_i(X)\times N^i(Z)_{\R}\mapsto (\alpha\cdot \beta)\in \R$$
is a continuous bilinear form.

\subsection{Positive generated cycles}
Define $\sG^+_i(X)$ to be the closure of $\Eff_i(X)$ in $\sG_i(X)$. It is a closed convex cone of $\sG_i(X).$ We call the induced topology on $\sG^+(X)$ the \emph{weak topology}.
We view $\sG_i^+(X)$ as an analogy of the notion of positive closed currents of bidimension $(i,i)$ in complex geometry.
It is clear that the projection $R_{X}: \sG_i(X)\to N_i(X)_{\R}$ maps $\sG^+_i(X)$ onto $\Psef_i(X)$.
Set $$R^+:=R_{X}|_{\sG^+_i(X)}: \sG^+_i(X)\twoheadrightarrow \Psef_i(X).$$
\begin{lem}\label{lempropcon}The map $R^+: \sG^+_i(X)\twoheadrightarrow \Psef_i(X)$ is proper. In other words, for every ample line bundle $L$ of $X$ and $B\geq 0,$
$$\{\alpha\in \sG^+_i(X)|\,\, (\alpha\cdot L^i)\leq B\}$$ is compact. In particular, $(R^+)^{-1}(0)=0.$
\end{lem}
\proof
For every $V\in \sZ(X),$
set $K_V:=\{\beta\in \Psef_i(V)|\,\, (\beta\cdot L^i)\leq B\}$, which is compact. 
Then we have $$\{\alpha\in \sG^+_i(X)|\,\, (\alpha\cdot L^i)\leq B\}\subseteq \prod_{V\in \sZ(X)}K_V.$$
By Tychonoff theorem, the right hand side is compact. As  $\{\alpha\in \sG^+_i(X)|\,\, (\alpha\cdot L^i)\leq B\}$ is close, it is compact.
As $$\R_{>0}\times (R^+)^{-1}(0)=(R^+)^{-1}(0)$$ and $(R^+)^{-1}(0)$ is compact, we get $(R^+)^{-1}(0)=0.$
This concludes our proof. 
\endproof

For every constructible subset $V$ of $X$, define $\sG^+_i(V,X):=\sG_i(V,X)\cap \sG_i^+(X).$
We may check that $\pi_V(\Eff_i(X))\subseteq \Eff_i(V,X).$ By the continuity, we get 
\begin{equation}\label{projpposi}\pi_V(\sG^+_i(X))=\sG^+_i(V,X).
\end{equation}
Indeed, $\sG^+_i(V,X)$ is the closure of $\Eff_i(V,X)$ in $\sG^+(X)$ and is exactly the space of $\alpha\in \sG^+_i(X)$ supported in $V.$
If $V$ is locally closed,  we may check that 
\begin{equation}\label{equationlocallose}R_{\overline{V}^{zar}}(\sG^+_i(V,X))=\Psef_i(V,X).
	\end{equation}

\medskip

We then have the following application on the positivity of some intersection numbers.
\begin{pro}\label{propositiveinter}Let $D$ be an effective Cartier divisor of $X$.
	Let $\alpha\in \sG^+_1(X)$. Assume that $\Supp\, \alpha\cap D=\emptyset,$ then $(\alpha\cdot D)\geq 0.$
\end{pro}
\proof
As $\alpha\in \sG^+_i(X\setminus D,X)$, we have $R_X(\alpha)\subseteq \Psef_i(X\setminus D,X).$ Hence we have $(\alpha\cdot D)\geq 0.$
\endproof

\subsubsection{Strong topology}
Define $\sD\sG_i^+:=\sG_i^+-\sG_i^+=\{\alpha-\beta|\,\, \alpha,\beta\in \sG_i^+\}$ which is a subspace of $\sG_i(X)$.
Let $L$ be any ample line bundle on $X.$ For every $\alpha\in \sD\sG_i^+,$ define 
$$\|\alpha\|_{L}:=\inf\{((\alpha_1+\alpha_2)\cdot L^i)|\,\,\alpha_1,\alpha_2\in \sG_i^+, \alpha_1-\alpha_2=\alpha\}.$$
Easy to check that $\|\cdot\|$ is a norm. Moreover, for any different ample line bundles $L_1,L_2$, the norms $\|\cdot\|_{L_1}, \|\cdot\|_{L_2}$ are equivalent i.e. 
there is $C>1$ such that 
$$C^{-1}\|\cdot\|_{L_2}\leq \|\cdot\|_{L_1}\leq C\|\cdot\|_{L_2}.$$
For the simplicity, we now fix $L$ and write $\|\cdot\|$ for $\|\cdot\|_{L_1}.$
Easy to check that for every $\alpha\in \sG_i^+$, we have 
\begin{equation}\label{equationconemetric}\|\alpha\|=(\alpha\cdot L^i).
\end{equation}

\medskip
We call the topology on $\sG^+_i(X)$ induced by $\|\cdot\|$ the strong topology.
As the map $\alpha\mapsto (\alpha\cdot L^i)=\|\alpha\|$ is continuous on $\sG_i^+(X)$, the map $$\id: (\sD\sG_i^+,\|\cdot\|)\to \sG^+_i(X)$$ is continuous.
So the strong topology on $\sG^+(X)$ is stronger than the weak topology.

The following result shows that $(\sD\sG_i^+,\|\cdot\|)$ is a Banach space. 
\begin{pro}\label{proconvergecone}Let $\alpha_n\in \sG_i^+$. The followings are equivalent:
\begin{points}
\item $\sum_{n\geq 0}\|\alpha_n\|<+\infty$;
\item $\sum_{n\geq 0}\alpha_n$ converges for weak topology;
\item $\sum_{n\geq 0}\alpha_n$ converges for $\|\cdot\|.$
\end{points}
In particular,  the norm $\|\cdot\|$ is complete on $\sD\sG_i^+$.
\end{pro}
\proof
It is clear that (ii) implies (i) and (iii) implies (i).

Assume (i), set  $B:=\sum_{n\geq 0}\|\alpha_n\|$. By Lemma \ref{lempropcon}, $K:=\{\alpha\in \sG^+_i(X)|\,\, (\alpha\cdot L^i)\leq B\}$ is compact. For $m\geq 0$, set $S_m:=\sum_{n=0}^m \alpha_n.$
The $S_m\in K$ for every $m\geq 0.$
Pick $\beta$ in the limit set. of $\{S_m|\,\, m\geq 0\}.$  We first show that $\sum_{n\geq 0}\alpha_n\to \beta$ in weak topology.  For every $m\geq 0$,  $\beta-S_m$ is contained in the limit set of $$\{S_{m+n}-S_m|\,\, n\geq 0\}\subseteq K\subseteq \sR^+_i(X).$$
Then we have $$\beta-S_m\in K$$ for every $m\geq 0.$
For every $V\in \sZ(X)$, we have $R_V(\beta-S_m)\subseteq \Psef_i(V).$
As $$0\leq ((\beta-S_m)\cdot L^i)\leq \sum_{n\geq m}\|\alpha_n\|\to 0,$$
we get 
$R_V(\beta-S_m)\to 0$ as $m\to \infty.$ So we have $\sum_{n\geq 0}\alpha_n\to \beta$ as $n\to \infty.$ Then we get that (i) implies (iii).
As $\beta-S_m\in K$, $\|\beta-S_m\|=((\beta-S_m)\cdot L^i)\to 0.$ So $\sum_{n\geq 0}\alpha_n\to \beta$ in $(\sD\sG_i^+, \|\cdot\|).$
So (i) implies (iii).
This concludes the proof.
\endproof

\medskip

The following example shows that the stronger topology is strictly stronger than the weak topology in general.

\begin{exe}\label{exeststrtweak}
Let $X=\P^2_{\C}.$ Let $Z_n, n\geq 0$ be the line defined by $y-x-nz=0.$ Let $L:=\sO_{\P^2_{\C}}(1)$.
It is clear that $Z_n$ converges weakly. But for every distinct $n,m\geq 0$, $\|Z_n-Z_m\|=1$. So there is no convergence subsequence for the strong topology.
\end{exe}

However, as the morphism $\alpha\mapsto (\alpha\cdot L^i)=\|\alpha\|$ is continuous on $\sG_i^+(X)$, a sequence $\alpha_n\in \sG_i^+(X)$ tends to $0$ weakly  if and only if it tends to $0$ strongly.

\subsubsection{Induced measure}
Denote by ${\rm P}^i(X)$ the dual cone of $\Psef_i(X)$ in $N^i(X)_{\R}.$
As $\Psef_i(X)$ has non-empty interior and is salient, $\sP^i(X)$ has non-empty interior and is salient.
We have the following morphism 
$$\nu: {\rm P}^i(X)\times \sG^+_i(X)\to \sM(|X|)$$ as follows:
For every $\beta\in {\rm P}^i(X)$, $\alpha\in \sG^+_i(X)$ and every constructible set $V\in \sA(|X|)$, define 
$$\nu(\beta, \alpha)(V):=(\pi_V(\alpha)\cdot \beta)\geq 0.$$
We define $\nu(\beta,\alpha)$ to be the unique measure in $\sM(|X|)$ make the above equality holds for every $V\in \sA(|X|).$

Such measure exists and is unique. Indeed by Carath\'eodory's extension theorem, we only need to show the following:
Let $V_{n}, n\geq 0$ be a sequence of disjoint constructible sets. Let $V\in \sA(|X|)$ satisfying $V=\sqcup_{n\geq 0}V_n$.
Then $$\nu(\beta,\alpha)(V)=\sum_{n\geq 0}\nu(\beta,\alpha)(V_n).$$
We note that constructible sets are open and closed. The compactness of $V$ implies that $V_n=\emptyset$ for $n\gg 0$.
By (\ref{equationdisjdec}), we get $\nu(\beta,\alpha)(V)=\sum_{n\geq 0}\nu(\beta,\alpha)(V_n).$

\begin{lem}\label{leminducemeacon}The map $$\nu: {\rm P}^i(X)\times \sG^+_i(X)\to \sM(|X|)$$ is continuous. Here $\sM(|X|)$ is endowed with the weak-* topology and $\sG^+_i(X)$ is endowed with the weak topology.
	\end{lem}
\proof
For the continuity, we only need to show that for every continuous function $h$ on $|X|$, $\nu^*h$ is continuous on ${\rm P}^i(X)\times \sG^+_i(X).$ By Lemma \ref{lemdensesetcon}, we may assume that $h$ takes form $1_U$ where $U\in \sA(|X|)$. Then the function 
$\mu^*h$ sends $(\beta,\alpha)$ to $(\beta\cdot \pi_V(\alpha))$ which is continuous.
\endproof

Then every $\alpha\in \sG^+_i(X)$ defines a unique vector-valued measure $\nu_{\alpha}^X$ such that for every $\beta\in {\rm P}^i(X)$,
$$\nu(\beta\cdot \alpha)=(\beta\cdot \nu_{\alpha}^X).$$
Easy to check that $\Supp\, \alpha=\Supp\, \nu_{\alpha}^X.$

\subsubsection{Atomic decomposition}
For $x\in |X|$, define $$\Psef_i(x,X):=\cap_U \Psef^i(U,X)$$ in $N_i(Z_x)_{\R}$, where $U$ is taken over all non-empty Zariski open subsets of $Z_x.$ If $\dim Z_x\leq i-1$, it is clear that $\Psef_i(x,X)=\{0\}.$
\begin{lem}\label{lemgermpseud}For $x\in |X|$ with $\dim Z_x=l\geq i$, then $\Psef_i(x,X)$ is a closed, convex, salient cone with non-empty interior.
	\end{lem}
\proof
Pick $L$ an ample line bundle on $X$.
For every non-empty Zariski open subset $U$ of $Z_x$, set $K_U:=\{v\in \Psef^i(U,X)|\,\, (v\cdot L^i)\leq 1\}.$ These $K_U$ are compact, as $\Psef^i(U,X)\subseteq \Psef^i(Z_x)$ and $K_{Z_{x}}$ is compact. Hence $\{v\in \Psef^i(x,X)|\,\, (v\cdot L^i)\leq 1\}=\cup_{U}K_U$ is compact.
Hence $\Psef_i(x,X)$ is closed. By De Jong's alteration theorem \cite{Jong1996}, there is a smooth projective variety $Y$ with a generically finite and surjective morphism $q:=Y\twoheadrightarrow Z_x$. Let ${\rm BPF}^{l-i}(Y)$ be the cone in $N^{l-i}(Y)_{\R}=N_i(Y)_\R$ as in \cite[Definition 3.3.1]{Dang2020}. It definition shows that ${\rm BPF}^{l-i}(Y)\subseteq \Psef_i(\eta,Y)$ where $\eta$ is the generic point of $Y$. Hence we have $q_*({\rm BPF}^{l-i}(Y))\subseteq \Psef_i(x,X).$ By \cite[Theorem 3.3.3 (1)]{Dang2020}, ${\rm BPF}^{l-i}(Y)$ has non-empty interior. As $q_*: N_i(Y)_{\R}\to N_i(X)_{\R}$ is surjective, $q_*({\rm BPF}^{l-i}(Y))$ (hence $\Psef_i(x,X)$) has non-empty interior.
As $\Psef_i(x,X)\subseteq \Psef_i(X)$ and $\Psef_i(X)$ is salient, $\Psef_i(x,X)$ is salient.
\endproof

For every $v\in \Psef_i(x,X)$, we define a positive generated cycle $v\delta_x$ of $X$ as the element 
$(\alpha_Z)_{Z\in \sZ(X)}\in \prod_{Z\in \sZ(X)}N_i(Z)_{\R}$ such that for every $Z\in \sZ(X)$, $\alpha_Z=0$ if $x\not\in Z$ and $\alpha_Z=v$ if $x\in Z.$ 
We now check that $v\delta_x$ is contained in $\sG^+_i(X)$. Let $Z_1,\dots, Z_m\in \sZ(X)$ with $x\in\cup_{j=1}^mZ_j$. Let $J$ be the set of $j$ such that $x\not\in Z_j.$ Set $U:=Z_x\setminus(\cup_{j\in J}Z_j)$ which is a non-empty Zariski open subset of $Z_x$. There are effective $i$-cycles $W_n$ taking form $W_n=\sum_{s=1}^la_{n,s}Z_{w_{n,s}}$ such that $a_{n,s}>0$, $w_{n,s}\in U$ and $[W_n]\to v$ in $N_i(Z_x)_{\R}.$
Then the image of $v\delta_x$ in $\prod_{j=1}^mN_i(Z_j)_{\R}$ can be approximated by the images of $W_n, n\geq 0.$
So $v\delta_x\in \prod_{Z\in \sZ(X)}N_i(Z)_{\R}$ is contained in the closure of effective $i$-cycles of $X$. Hence we get $v\delta_x\in \sG^+_i(X)$.  we call such $v\delta_x$ the \emph{atoms} in $\sG^+_i(X)$

\medskip

The following result shows that every positive generated cycle is a positive combination of at most countably many atoms.
\begin{thm}\label{thmdecogc}Every $\alpha\in \sG^+_i(X)$ takes form 
	$$\alpha=\sum_{j\geq 0}v_j\delta_{x_j}$$
	where $v_j\in \Psef_i(x_j,X)$ with $\sum_{j\geq 0}(v_j\cdot L^i)<+\infty$ where $L$ is any ample line bundle on $X.$
\end{thm}

Combing Theorem \ref{thmdecogc} with Corollary \ref{corsmxsc}, we indeed showed that $\sG^+_0(X)=\sM(|X|).$

\proof
By Theorem \ref{thmRadon} and Corollary \ref{corsmxsc}, write $\nu(L^i,\alpha)=\sum_{j\geq 0}a_j\delta_{x_j}.$ Set $Z_j:=Z_{x_j}.$ 
For each $j$, $\alpha_j:=\pi_{Z_j}\alpha\in \sP^+_i(Z_{j})$ defines a vector-valued measure 
$\nu_{\alpha_j}^{Z_j}\in \sM(|Z_j|,\Psef_i(Z_j)).$ We may write $\nu_{\alpha_j}^{Z_j}=v_j\delta_{x_j}+\beta_j$ where $\epsilon_j\in \sM(|Z_j|,\Psef_i(Z_j))$ with $\beta_j(\{x_j\})=0.$ We have $(v_j\cdot L^i)=a_j.$

For every $\epsilon>0$, there is a non-empty Zariski open subset $U_{\epsilon}$ of $Z_j$ such that for every non-empty Zariski open subset $W$ of $U_{\epsilon},$ we have 
$$(L^i\cdot (\pi_W(\alpha)-v_j))=\int_{W}(L^i\cdot \beta)<\epsilon.$$ 
As $R_{Z_j}(\pi_W(\alpha))\subseteq N_i(W,X)$, we get $v_j\in \Psef_i^+(x_j,X).$ 

\medskip
The above construction shows that for every $Z\in \sZ(X)$,
$$\nu^Z_{\pi_Z(\alpha)}=\nu^Z_{\pi_Z(\sum_{j\geq 0}v_j\delta_{x_j})}.$$
Then we have $$R_Z(\alpha)=\nu^Z_{\pi_Z(\alpha)}(Z)=R_Z(\sum_{j\geq 0}v_j\delta_{x_j}).$$
This concludes the proof.
\endproof

\begin{cor}\label{corsequecom}Let $L$ be an ample line bundle on $X$. Then $K:=\{\alpha\in \sG_i(X)|\,\, (L^i\cdot \alpha)\leq 1\}$ is compact and sequentially compact. 
	\end{cor}

\proof
By Lemma \ref{lempropcon}, $K$ is compact. We now show that $K$ is sequentially compact. Let $\alpha_n\in K, n\geq 0.$ Then we have $\nu(L^i,\alpha_n)\in \sM^1(|X|).$ By Corollary \ref{corsmxsc}, up to taking subsequences, we may assume that $\nu(L^i,\alpha_n)$ converges to a measure $$\nu=\sum_{j\geq 0}a_j\delta_{x_j}\in \sM^1(|X|).$$  We may assume the above $x_j$ are distinct. Set $Z_j:=Z_{x_j}.$
 By diagonal method, after taking subsequences, we may assume that $\nu^{Z_{j}}_{\pi_{Z_j}\alpha_n},n\geq 0$ converges.

We want to show that $\alpha_n$ converges. 
For this, we only need to show that for every $Z\in \sZ(X)$, $R_Z(\alpha_n)$ converges. For every $\epsilon>0$, there is $M\geq 0$ such that $\sum_{j\geq M}a_j<\epsilon/2.$ For every $j=1,\dots, M$, pick a non-empty Zariski open subset $U_j$ of $Z_j$ such that for every $s\neq j$ with $x_j\not\in Z_s$, we have $U_j\cap Z_s=\emptyset.$ Then $U_1,\dots, U_M$ are disjoint constructible subsets.  After shrinking $U_j$, we may assume that $U_j\cap Z=\emptyset$ if and only if $x_j\not\in Z.$ Set $U_0:=X\setminus (\cup_{j=1}^MU_j).$ 
 We have $$R_Z(\alpha_n)=\sum_{j=1}^M \nu^Z_{\pi_{Z_j}(\alpha_n)}(U_j)+\nu^Z_{\alpha_n}(U_0).$$
For every $j=1,\dots,M$,  $\nu^Z_{\pi_{Z_j}(\alpha_n)}$ converges, so there is $N\geq 0$ such that for $n\geq m\geq N$, we have
$\|\nu^Z_{\pi_{Z_j}(\alpha_n)}(U_j)-\nu^Z_{\pi_{Z_j}(\alpha_m)}(U_j)\|_L\leq \epsilon/2M$, where $\|\cdot\|_L$ is the norm on $N_i(Z)_{\R}$ induced by $L^{i}\in {\rm P}^i(Z)^{\circ}.$ Moreover $\|\nu^Z_{\alpha_n}(U_0)\|_L\leq \sum_{j\geq M}a_j\leq \epsilon/2.$
So we get $$\|R_Z(\alpha_n)-R_Z(\alpha_m)\|_L\leq \epsilon$$ for $n\geq m\geq N.$ This concludes the proof.
\endproof

\subsection{A dynamical application}
Let $f: X\to X$ be a surjective endomorphism. Let $L$ be an ample line bundle on $X$.

\begin{pro}\label{prowanderingclass}Let $x\in X$ and $v\in \Psef_1(x,X)$. Set $\alpha:=v\delta_x$. Let $D$ be any effective Cartier divisor on $X$. Assume that the orbit $O_f(x)$ is Zariski dense in $X.$ Then we have 
	\begin{equation}\label{equationcomeffam}\liminf_{n\to \infty}\frac{(f^n_*(\alpha)\cdot D)}{(L\cdot f^n_*(\alpha))}\geq 0.
		\end{equation}
	
Moreover, if $M$ is a big line bundle, then there is $\delta>0$, such that $$(f^n_*(\alpha)\cdot M)> \delta (f^n_*(\alpha)\cdot L)$$ for $n\gg 0.$
\end{pro}

\proof
There is $\beta>0$ and $s\geq 0$ such that $$(f^n_*(\alpha)\cdot D)=(\sum_{i=1}^mc_ie^{i\theta_in})\beta^nn^s+O(\beta^nn^{s-1})$$
where $c_i\neq 0$ and $\theta_i\in \R/2\pi\Z$ are distinct. 
It is clear that 
\begin{equation}\label{comparebetanorm}\beta^nn^s\leq C(f^n_*(\alpha)\cdot L)
	\end{equation}
	 for some $C>0.$
	 \medskip

Define a continuous function $h: (\R/2\pi\Z)^m\to \R$ by $$(\phi_1,\dots,\phi_m)\mapsto \sum_{i=1}^mc_ie^{i\phi_i}.$$
Define $\Theta: \Z\to (\R/2\pi\Z)^m$ by $$n\mapsto (\theta_1n,\dots,\theta_mn).$$
Let $Z$ the the closure of $\Theta(\Z).$ 
By Poincar\'e recurrence theorem,  there is a strictly increasing sequence $n_i\in \Z_{\geq 0}, i\geq 0$ such that $$\Theta(n_i)\to 0.$$ Then for every $m>0$, $$\Theta(m)=\lim_{i\to \infty}\Theta(m-n_i)$$
and
$$\Theta(-m)=\lim_{i\to \infty}\Theta(-m+n_i).$$
So $Z$ is also the closure of $\Theta(\Z_{\geq 0})$ and of $\Theta(\Z_{\leq 0}).$
\begin{lem}\label{lemfinitestu}For every non-empty open subset $U$ of $Z$, there is $r_U\geq 1$ such that for every $n\geq 0$, there is $m\in \{n,n+1,\dots, n+r_U\}$ such that $\Phi(m)\in U.$
	\end{lem}

\medskip
We claim that $h\geq 0$ on $Z$. Otherwise, there is $b<0$ such that $$U:=\{z\in Z|\,\, h(z)<b\}\neq \emptyset.$$
Set $$W:=\{n\geq 0|\,\, \Phi(n)\in U\}$$ and 
$$w(n):=\#(\{0,\dots,n-1\}\cap W).$$
By Lemma \ref{lemfinitestu}, we have 
\begin{equation}\label{equwnongzero}\liminf_{n\to \infty}w(n)/n>0.
\end{equation}
As
\begin{equation}\label{equationfnaind}(f^n_*(\alpha)\cdot D)=h(\Phi(n))\beta^nn^s+O(\beta^nn^{s-1}),
	\end{equation}
for $n\in W$, we have $$(f^n_*(\alpha)\cdot D)<0.$$
Then, by Proposition \ref{propositiveinter}, we get $f^n(x)\in \Supp\, D$ for every $n\in W.$
As $O_f(x)$ is Zariski dense in $X$, by the Weak dynamical Mordell-Lang \cite[Theorem 1.17]{Xie2023} (see also \cite[Theorem 2.5.8]{Favre2000a}), we have 
$$\lim\#\{m=0,\dots,n-1|\,\, f^n(x)\in \Supp\, D\}/n=0.$$
It follows that $$\lim_{n\to \infty}w_n/n=0,$$ which contradicts to (\ref{equwnongzero}).
This proves the claim.

\medskip

By (\ref{equationfnaind}), we get $$(f^n_*(\alpha)\cdot D)\geq -C'\beta^nn^{s-1}$$ for some $C'>0.$
We conclude the proof of (\ref{equationcomeffam}) by (\ref{comparebetanorm}).

\medskip

As $M$ is big, we may write $M=\delta' L+E$ where $\delta'>0$ and $E$ is effective.
By (\ref{equationcomeffam}), for $n\gg 0$, we have $$(f^n_*(\alpha)\cdot E)>-\frac{\delta'}{2}(f^n_*(\alpha)\cdot L)$$
Hence we have  $$(f^n_*(\alpha)\cdot M)>\delta'(f^n_*(\alpha)\cdot L)-\frac{\delta'}{2}(f^n_*(\alpha)\cdot L)=\frac{\delta'}{2}(f^n_*(\alpha)\cdot L).$$
This concludes the proof.
\endproof

\proof[Proof of Lemma \ref{lemfinitestu}]
There is $r\geq 0$ such that $\Phi(r)\in U.$ We only need to Lemma \ref{lemfinitestu} for the open subset $U-\Phi(r)$. After replacing $U$ by $U-\Phi(r)$, we may assume that $0\in U.$ After replacing $U$ by $U\cap (-U)$, we may assume that $U=-U.$ Pick an open neighborhood $V$ of $0$ such that $V=-V$ and $V+V\subseteq U.$ As $Z$ is compact, there is a finite subset $F$ of $Z$ such that $F+V=Z.$ For every $z\in F$, there is $r_z\geq 0$ such that $\Phi(-r_z)\in z+V.$ In other words, we have $$z+\Phi(r_z)\in -V=V.$$
For every $n\geq 0$, there is $z_n\in F$ such that $\Phi(n)\in z_n+V.$ Then we have $$\Phi(n+r_{z_n})\in z_n+V+\Phi(r_{z_n})\in V+V\subseteq U.$$ Set $r_U:=\max\{r_{z}|\,\, z\in F\}.$
We conclude the proof. 
\endproof

\section{The spectrum for the ample cone}
Let $X$ be a projective variety over $\bk$ of dimension $d$ and $f: X\to X$ is a surjective endomorphism. 
Let $V$ be an irreducible and periodic subvariety of $X$ of dimension $d_V\geq 0$. Recall that $$\mu_i(V,f):=\mu_i(f|_V^{r_V}:V\to V)^{1/l}$$ where $r_V\geq 1$ is a period of $V$. It does not depend on the choice of $l.$

\medskip

%
%
%

The aim of this section is to prove the following result.
\begin{thm}(=Theorem \ref{thmquaisamptoamintr})\label{thmquaisamptoam}
	We have $$\Sp(f^*,\Amp(X))=\cup_{V}\{\mu_i(V,f) |\,\, i=1,\dots, d_V\}$$ where the union taken over all irreducible periodic subvarieties.

	In other words,	for $\alpha\in \R_{>0}$, $f$ is $\alpha$-amplified if and only if for every periodic irreducible subvariety $V$, $f^{r_V}|_V$ is $\alpha$-amplified, where $r_V\geq 1$ is a period of $V$. 
\end{thm}

\subsection{Growth rate}
Let $\|\cdot\|$ be any norm on $N_1(X)_{\R}.$ For every $v\in N_1(X)_{\R}\setminus \{0\}$, by (\ref{equnormgrowth})  in Section \ref{sectionlinearalg}, $\beta_f(v):=\lim\limits_{n\to\infty}\|f^n_*(v)\|^{1/n}$ converges and does not depend on the choice of $\|\cdot\|$.
We call it the \emph{growth rate} of $v.$
By (\ref{equnspgrothrate}) in Section \ref{sectionlinearalg}, $\beta_f(v)\in \{|c||\,\, c\in \Sp(f)\}.$ So $\beta_f(v)$ has only finitely many possible values. Moreover $\beta_f(\alpha)\leq \la_1(f)=\mu_1(f).$

\medskip

The aim of this section is to study $\beta_f(v)$ for $v\in \Psef_1(X)\setminus \{0\}.$ 
Let $L$ be any ample line bundle on $X$. By Remark \ref{remconenorm}, $L$ induces a norm $\|\cdot\|_L$ on $N^1(X)_{\R}$ such that for every $v\in \Psef_1(X)$, $\|v\|_L=(v\cdot L).$ Using this norm, for $v\in \Psef_1(X)$, we get $$\beta_f(v)=\lim_{n\to\infty}(f^n_*(v)\cdot L)^{1/n}.$$

\medskip

\begin{lem}\label{lemgrforpsedadd}Let $v_i\in \Psef_1(X), i\geq 0$ with $\sum_{i\geq 0}\|v_i\|<+\infty$. Set $v:=\sum_{i\geq 0}v_i$. Assume that $v\neq 0$, then 
	we have 
	$$\beta_f(v)=\max\{\beta_f(v_i)|\,\, v_i\neq 0\}.$$
\end{lem}
\proof
For every $i\geq 0$ with $v_i\neq 0$, we have $$(f_*^n(v)\cdot L)\geq (f_*^n(v_i)\cdot L).$$
Hence $\rho_f(v)\geq \rho_f(v_i).$ As $\rho_f(v_i)$ has only finitely many possible values, we get 
$$\beta_f(v)\geq \max\{\beta_f(v_i)|\,\, v_i\neq 0\}.$$
Set $\beta:=\max\{\beta_f(v_i)|\,\, v_i\neq 0\}.$ 
By (\ref{equnspgrothrate}) in Section \ref{sectionlinearalg}, we have $v_i\in E_{\overline{\D}(\beta)}$ for every $i\geq 0$. As $E_{\overline{\D}(\beta)}$ is closed, we have $v\in E_{\overline{\D}(\beta)}.$ Hence $\rho_f(v)\leq \beta.$
This concludes the proof. 
\endproof

Every $v\in \Psef_1(X)$ can be presented by an positive generated cycle $\alpha\in \sG^+_1(X).$ 
For the simplicity, we also write $\rho_f(\alpha)$ for $\rho_f(v).$
By Theorem \ref{thmdecogc}, $\alpha$ is a positive combination of at most countably many atoms.
By Lemma \ref{lemgrforpsedadd}, we only need to understand the growth rate of atoms in $\sG^+_1(X).$
\begin{pro}\label{proatomgrowth}Let $x\in X$ and $v\in \Psef_1(x,X)$. Set $\alpha:=v\delta_x$.  Assume that the orbit $O_f(x)$ is Zariski dense in $X.$ Then we have $$\beta_f(\alpha)\in \{\mu_i(f), i=1,\dots, d\}.$$
\end{pro}
\proof
Assume by contradiction that $\beta_f(\alpha)\not\in  \{\mu_i(f), i=1,\dots, d\}$. Set $\mu_{d+1}=0.$ There is a unique $i=1,\dots, d$ such that  $\beta_f(\alpha)\in (\mu_i, \mu_{i+1}).$

\medskip

For every $n\geq 0$, define $L_n:=(f^n)^*L$ and $\alpha_n:=(f^n)_*\alpha.$ By projection formula,  we have 
$$(L_{n_1}\cdot \alpha_{n_2})=(L\cdot \alpha_{n_1+n_2}).$$
Pick $\epsilon\in (0,1)$ such that $$\epsilon^{-1}\mu_{i+1}<\beta< \epsilon^2\mu_i.$$
There is $m_0\geq 1$ such that for every $m\geq m_0$,
$$\epsilon^{2m}\mu_i^m+\epsilon^{-m}\mu_{i+1}^m<\epsilon^m\mu_i^m.$$
By \cite[Theorem 3.7]{Xie2024},  there is $m\geq 1$ such that
$$M:=L_{2m}+\mu_i^{m}\mu_{i+1}^mL-\epsilon^m\mu_i^{m}L_m$$ is big.
By Proposition \ref{prowanderingclass}, there is $N\geq 0$ such that for every $n\geq N$,
$$(\alpha_{mn}\cdot M)\geq 0.$$
Then we get 
$$(\alpha_{(n+2)m}\cdot L)+\mu_i^{m}\mu_{i+1}^m(\alpha_{nm}\cdot L)-\epsilon^m\mu_i^{m}(\alpha_{(n+1)m}\cdot L)\geq 0$$
for $n\geq N.$
It follows that 
\begin{equation}\label{equrec}(\alpha_{(n+2)m}\cdot L)-\epsilon^{-m}\mu_{i+1}^m(\alpha_{(n+1)m}\cdot L)\geq \epsilon^{2m}\mu_i^m((\alpha_{(n+1)m}\cdot L)-\epsilon^{-m}\mu_{i+1}^m(\alpha_{nm}\cdot L))
\end{equation}
for all $n\geq N.$
As $\beta_f(\alpha)=\beta>\epsilon^{-1}\mu_{i+1}$, we have 
$$(\alpha_{(n+1)m}\cdot L)-\epsilon^{-m}\mu_{i+1}^m(\alpha_{nm}\cdot L)>0$$ for $n\gg 0.$
Then we get $$\beta^m=\liminf_{n\to \infty} (\alpha_{nm}\cdot L)^{1/n}\geq \epsilon^{-m}\mu_{i+1}^m,$$
which is a contradiction.
\endproof
Apply Proposition \ref{proatomgrowth} for every periodic irreducible subvarieties, we get the following result for any atom.
\begin{cor}\label{corgeneralatom}
	Let $x\in X$ and $v\in \Psef_1(x,X)$. Set $\alpha:=v\delta_x$. Then every irreducible component $V$  of $\cap_{m\geq 0}\overline{O_f(f^m(x))}^{zar}$ is $f$-periodic.
	Then we have $$\beta_f(\alpha)\in \{\mu_i(V,f), i=1,\dots, d\}.$$
	\end{cor}

Combine Lemma \ref{lemgrforpsedadd}, Corollary \ref{corgeneralatom} and Theorem \ref{thmdecogc}, we get the following result.
\begin{thm}\label{thmpsedgrowrate} For every $v\in \Psef_1(X)\setminus \{0\}$, we have 
	$$\beta_f(v)\in \cup_{V}\{\mu_i(V,f) |\,\, i=1,\dots, d_V\}$$
	where the union taken over all irreducible periodic subvarieties.
	\end{thm}

\subsection{The spectrum for the ample cone}

\proof[Proof of Theorem \ref{thmquaisamptoam}]
For $\alpha\in \R_{>0}$, if $f$ is $\alpha$-amplified, then for every periodic irreducible subvariety $V$ with period $r_V$, $f^{r_V}|_{V}$ is $\alpha^{r_V}$-amplified. 
Hence we get 
$$\cup_{V}\Sp((f^{r_V}|_V)^*, \Amp(V))^{1/r_V}\subseteq \Sp(f^*,\Amp(X)),$$
here $\Sp((f^{r_V}|_V)^*, \Amp(V))^{1/r_V}:=\{\beta^{1/r_V}|\,\, \beta\in \Sp(f^{r_V}|_V, \Amp(V))\}.$
As the big cone contains the ample cone, we have
$$\{\mu_i(V,f) |\,\, i=1,\dots, d_V\}\subseteq \Sp((f^{r_V}|_V)^*, \Amp(V))^{1/r_V}.$$
So we get 
$$\cup_{V}\{\mu_i(V,f) |\,\, i=1,\dots, d_V\}\subseteq \Sp(f^*,\Amp(X)).$$

\medskip

Set $S:=\cup_{V}\{\mu_i(V,f) |\,\, i=1,\dots, d_V\}.$
By contradiction, assume that $$\Sp(f^*,\Amp(X))\not\subseteq S.$$
By Theorem \ref{thmceighear}, $E_S\cap \Amp(X)=\emptyset.$ 
By Hahn-Banach theorem, there is $Z\in \Psef_1(X)\setminus \{0\}$ such that 
 $E_S\subseteq Z^{\bot}$. In particular, we have $\beta_f(Z)\not\in S.$ This contradicts to 
Theorem \ref{thmpsedgrowrate}.
\endproof

\begin{cor}\label{corfactoramplified}Let $Y$ be a projective variety over $\bk$ and $g:Y\to Y$ be an endomorphism. Let $\pi: X\to Y$ be a surjective morphism such that $\pi\circ f=g\circ \pi.$ If $f$ is $\alpha$-amplified for some $\alpha\in \R_{>0}$, then $g$ is $\alpha$-amplified.
	\end{cor}
\proof
Assume that $f$ is $\alpha$-amplified.
By contradiction, assume that $g$ is not $\alpha$-amplified. By Theorem \ref{thmquaisamptoam}, there is an irreducible $g$-periodic subvariety $V$ of $Y$ such that $\alpha\in \{\mu_i(V,g)|\,\, i=1,\dots, d_V\}.$ After replacing $f,g$ by a suitable iterate, we may assume that $g(V)=V.$
There is an irreducible component $W$ of $\pi^{-1}(V)$ which is $f$-periodic and satisfies $\pi(W)=V$.  =After replacing $f,g$ by a suitable iterate, we may assume that $f(W)=W.$
The product formula for relative dynamical degrees (c.f. \cite{Dinh2011}, \cite{Dang2020} and \cite[Theorem 1.3]{Truong2020}) shows that $\{\mu_i(V,g)|\,\, i=1,\dots, d_V\}\subseteq \{\mu_i(W,f)|\,\, i=1,\dots, d_W\}.$ Theorem \ref{thmquaisamptoam}, we have 
$$\alpha\in \{\mu_i(W,f)|\,\, i=1,\dots, d_W\}\subseteq \Sp(f^*,\Amp(X)).$$ The $f$ is not $\alpha$-amplified, which is a contradiction.
\endproof

	\bibliography{dd}
\end{document}